\documentclass[a4paper,twoside]{article}
\usepackage{a4}
\usepackage{amssymb}
\usepackage{amsmath}
\usepackage{upref}
\usepackage{url}
\usepackage[active]{srcltx}
\usepackage[pagebackref,colorlinks,citecolor=blue,linkcolor=blue,urlcolor=blue]{hyperref}
\usepackage[dvipsnames]{color}
\allowdisplaybreaks[2] 
\newcount\minutes \newcount\hours
\hours=\time
\divide\hours 60
\minutes=\hours
\multiply\minutes -60
\advance\minutes \time
\newcommand{\klockan}{\the\hours:{\ifnum\minutes<10 0\fi}\the\minutes}
\newcommand{\tid}{\today\ \klockan}
\newcommand{\prtid}{\smash{\raise 10mm \hbox{\LaTeX ed \tid}}}
\renewcommand{\prtid}{}
\makeatletter
\def\sectionmark#1{} 
\def\subsectionmark#1{}
\newcommand{\sectnr}{\ifnum \c@secnumdepth >\z@
                 \thesection.\hskip 1em\relax \fi}
\def\@evenhead{\footnotesize\rm\thepage\hfil\leftmark\hfil\llap{\prtid}}
\def\@oddhead{\footnotesize\rm\rlap{\prtid}\hfil\rightmark\hfil\thepage}
\def\tableofcontents{\section*{Contents} 
 \@starttoc{toc}}
\makeatother
\makeatletter
\def\@biblabel#1{#1.}
\makeatother
\makeatletter
\let\Thebibliography=\thebibliography
\renewcommand{\thebibliography}[1]{\def\@mkboth##1##2{}\Thebibliography{#1}
\addcontentsline{toc}{section}{References}
\frenchspacing 
\setlength{\@topsep}{0pt}
\setlength{\itemsep}{0pt}%
\setlength{\parskip}{0pt plus 2pt}%
}
\makeatother
\makeatletter
\def\mdots@{\mathinner.\nonscript\!.%
 \ifx\next,.\else\ifx\next;.\else\ifx\next..\else
 \nonscript\!\mathinner.\fi\fi\fi}
\let\ldots\mdots@
\let\cdots\mdots@
\let\dotso\mdots@
\let\dotsb\mdots@
\let\dotsm\mdots@
\let\dotsc\mdots@
\def\vdots{\vbox{\baselineskip2.8\p@ \lineskiplimit\z@
    \kern6\p@\hbox{.}\hbox{.}\hbox{.}\kern3\p@}}
\def\ddots{\mathinner{\mkern1mu\raise8.6\p@\vbox{\kern7\p@\hbox{.}}%
    \raise5.8\p@\hbox{.}\raise3\p@\hbox{.}\mkern1mu}}
\makeatother
\makeatletter
\let\Enumerate=\enumerate
\renewcommand{\enumerate}{\Enumerate%
\setlength{\@topsep}{0pt}
\setlength{\itemsep}{0pt}%
\setlength{\parskip}{0pt plus 1pt}%
\renewcommand{\theenumi}{\textup{(\alph{enumi})}}%
\renewcommand{\labelenumi}{\theenumi}%
}
\let\endEnumerate=\endenumerate
\renewcommand{\endenumerate}{\endEnumerate\unskip}
\makeatother
\makeatletter
\def\@seccntformat#1{\csname the#1\endcsname.\quad}
\makeatother
\newcommand{\authortitle}[2]{\author{#1}\title{#2}\markboth{#1}{#2}}
\newcommand{\auth}[2]{{#1, #2.}}
\newcommand{\art}[6]{{\sc #1, \rm #2, \it #3 \bf #4 \rm (#5), \mbox{#6}.}}

\newcommand{\artprep}[3]{{\sc #1, \rm #2, #3.}}
\newcommand{\arttoappear}[3]{{\sc #1, \rm #2, to appear in \it #3}}
\newcommand{\book}[3]{{\sc #1, \it #2, \rm #3.}}
\newcommand{\AND}{{\rm and }}
\newcommand{\artnopt}[6]{{\sc #1, \rm #2, \it #3\/ \bf #4 \rm (#5), \mbox{#6}}}
\RequirePackage{amsthm}
\newtheoremstyle{descriptive}%
  {\topsep}   
  {\topsep}   
  {\rmfamily} 
  {}          
  {\bfseries} 
  {.}         
  { }         
  {}          
\newtheoremstyle{propositional}%
  {\topsep}   
  {\topsep}   
  {\itshape}  
  {}          
  {\bfseries} 
  {.}         
  { }         
  {}          
\newtheoremstyle{remarkstyle}%
  {\topsep}   
  {\topsep}   
  {\rmfamily}  
  {}          
  {\itshape} 
  {.}         
  { }         
  {}          
\theoremstyle{propositional}
\newtheorem{thm}{Theorem}[section]
\newtheorem{prop}[thm]{Proposition}
\newtheorem{lem}[thm]{Lemma}
\newtheorem{cor}[thm]{Corollary}
\theoremstyle{descriptive}
\newtheorem{deff}[thm]{Definition}
\newtheorem{remark}[thm]{Remark}
\makeatletter
\renewenvironment{proof}[1][\proofname]{\par
  \pushQED{\qed}%
  \normalfont
  \trivlist
  \item[\hskip\labelsep
        \itshape
    #1\@addpunct{.}]\ignorespaces
}{%
  \popQED\endtrivlist\@endpefalse
}
\makeatother
\newcommand{\setm}{\setminus}
\renewcommand{\subsetneq}{\varsubsetneq}
\renewcommand{\emptyset}{\varnothing}
\def\vint{\mathop{\mathchoice%
          {\setbox0\hbox{$\displaystyle\intop$}\kern 0.22\wd0%
           \vcenter{\hrule width 0.6\wd0}\kern -0.82\wd0}%
          {\setbox0\hbox{$\textstyle\intop$}\kern 0.2\wd0%
           \vcenter{\hrule width 0.6\wd0}\kern -0.8\wd0}%
          {\setbox0\hbox{$\scriptstyle\intop$}\kern 0.2\wd0%
           \vcenter{\hrule width 0.6\wd0}\kern -0.8\wd0}%
          {\setbox0\hbox{$\scriptscriptstyle\intop$}\kern 0.2\wd0%
           \vcenter{\hrule width 0.6\wd0}\kern -0.8\wd0}}%
          \mathopen{}\int}
{\catcode`p =12 \catcode`t =12 \gdef\eeaa#1pt{#1}}      
\def\accentadjtext#1{\setbox0\hbox{$#1$}\kern   
                \expandafter\eeaa\the\fontdimen1\textfont1 \ht0 }
\def\accentadjscript#1{\setbox0\hbox{$#1$}\kern 
                \expandafter\eeaa\the\fontdimen1\scriptfont1 \ht0 }
\def\accentadjscriptscript#1{\setbox0\hbox{$#1$}\kern   
                \expandafter\eeaa\the\fontdimen1\scriptscriptfont1 \ht0 }
\def\accentadjtextback#1{\setbox0\hbox{$#1$}\kern       
                -\expandafter\eeaa\the\fontdimen1\textfont1 \ht0 }
\def\accentadjscriptback#1{\setbox0\hbox{$#1$}\kern     
                -\expandafter\eeaa\the\fontdimen1\scriptfont1 \ht0 }
\def\accentadjscriptscriptback#1{\setbox0\hbox{$#1$}\kern 
                -\expandafter\eeaa\the\fontdimen1\scriptscriptfont1 \ht0 }
\def\itoverline#1{{\mathsurround0pt\mathchoice
        {\rlap{$\accentadjtext{\displaystyle #1}
                \accentadjtext{\vrule height1.593pt}
                \overline{\phantom{\displaystyle #1}
                \accentadjtextback{\displaystyle #1}}$}{#1}}
        {\rlap{$\accentadjtext{\textstyle #1}
                \accentadjtext{\vrule height1.593pt}
                \overline{\phantom{\textstyle #1}
                \accentadjtextback{\textstyle #1}}$}{#1}}
        {\rlap{$\accentadjscript{\scriptstyle #1}
                \accentadjscript{\vrule height1.593pt}
                \overline{\phantom{\scriptstyle #1}
                \accentadjscriptback{\scriptstyle #1}}$}{#1}}
        {\rlap{$\accentadjscriptscript{\scriptscriptstyle #1}
                \accentadjscriptscript{\vrule height1.593pt}
                \overline{\phantom{\scriptscriptstyle #1}
                \accentadjscriptscriptback{\scriptscriptstyle #1}}$}{#1}}}}
\def\itunderline#1{{\mathsurround0pt\mathchoice
        {\rlap{$\underline{\phantom{\displaystyle #1}
                \accentadjtextback{\displaystyle #1}}$}{#1}}
        {\rlap{$\underline{\phantom{\textstyle #1}
                \accentadjtextback{\textstyle #1}}$}{#1}}
        {\rlap{$\underline{\phantom{\scriptstyle #1}
                \accentadjscriptback{\scriptstyle #1}}$}{#1}}
        {\rlap{$\underline{\phantom{\scriptscriptstyle #1}
                \accentadjscriptscriptback{\scriptscriptstyle #1}}$}{#1}}}}
\newcommand{\Cpt}{{C_{\theta,p}}}
\newcommand{\CptY}{{C_{\theta,p}^{Y}}}
\newcommand{\cpt}{{\capp_{\theta,p}}}
\newcommand{\cptY}{{\capp_{\theta,p}^{Y}}}
\newcommand{\cptYprime}{{\capp_{\theta,p}^{Y'}}}
\newcommand{\cptYhat}{{\capp_{\theta,p}^{\Yhat}}}
\newcommand{\cptZ}{{\capp_{\theta,p}^{Z}}}
\newcommand{\cptX}{{\capp_{\theta,p}^{X}}}
\newcommand{\cptXeps}{{\capp_{\theta,p}^{\dXeps}}}
\newcommand{\CpXeps}{C_p^{\clXeps}}
\newcommand{\cpXeps}{\capp_p^{\clXeps}}
\newcommand{\BdXeps}{B^{\dXeps}}
\newcommand{\BXeps}{B^{\clXeps}}
\newcommand{\BZ}{B^Z}
\newcommand{\BY}{B}
\newcommand{\BYY}{B^Y}
\newcommand{\BYp}{B^{Y'}}
\newcommand{\BYhat}{B^{\Yhat}}
\DeclareMathOperator{\diam}{diam}
\DeclareMathOperator{\dist}{dist}
\DeclareMathOperator{\capp}{cap}
\newcommand{\cp}{\capp_p}
\newcommand{\bdry}{\partial}
\newcommand{\bdy}{\bdry}
\newcommand{\simge}{\gtrsim}
\newcommand{\simle}{\lesssim}
\newcommand{\al}{\alpha}
\newcommand{\alp}{\alpha}
\newcommand{\be}{\beta}
\newcommand{\ga}{\gamma}
\newcommand{\de}{\delta}
\newcommand{\eps}{\varepsilon}
\newcommand{\la}{\lambda}
\newcommand{\Om}{\Omega}
\newcommand{\p}{{$p\mspace{1mu}$}}
\newcommand{\R}{\mathbf{R}}
\def\cprime{{\mathsurround0pt$'$}}
\newcommand{\Np}{N^{1,p}}
\newcommand{\lQo}{\itunderline{Q}_0}
\newcommand{\lQoY}{\lQo^{Y}}
\newcommand{\uQo}{\itoverline{Q}_0}
\newcommand{\uQoY}{\uQo^{Y}}
\newcommand{\lSo}{\itunderline{S}_0}
\newcommand{\uSo}{\itoverline{S}_0}
\newcommand{\lQ}{\itunderline{Q}}
\newcommand{\uQ}{\itoverline{Q}}
\newcommand{\uq}{\overline{q}}
\newcommand{\lqq}{\underline{q}} 
\newcommand{\us}{\itoverline{s}}
\newcommand{\ls}{\itunderline{s}}
\newcommand{\uqo}{\uq_0}
\newcommand{\lqo}{\lqq_0}
\newcommand{\uso}{\us_0}
\newcommand{\lso}{\ls_0}
\newcommand{\Ga}{\Gamma}
\newcommand{\clX}{\itoverline{X}}
\newcommand{\clXeps}{\clX_\eps}
\newcommand{\mube}{{\mu_\be}}
\newcommand{\bdyeps}{\bdy_\eps}
\newcommand{\dXeps}{\bdyeps X}
\newcommand{\bdyXeps}{\dXeps}
\newcommand{\kard}{\hat{\ka}}
\newcommand{\kat}{\hat{\ka}}
\newcommand{\kaeps}{\ka_\eps}
\newcommand{\ka}{\kappa}
\DeclareMathOperator{\supp}{supp}
\DeclareMathOperator{\spt}{supp}
\newcommand{\lqoZ}{\lqq_0^Z}
\newcommand{\uqoY}{\uq_0^Y}
\newcommand{\uqoX}{\uq_0^{\clXeps}}
\newcommand{\lqoX}{\lqq_0^{\clXeps}}
\newcommand{\Yhat}{{\widehat{Y}}}
\newcommand{\dhat}{{\hat{d}}}
\newcommand{\nuhat}{{\hat{\nu}}}
\newcommand{\rt}{\tilde{r}}
\numberwithin{equation}{section}
\newcommand{\eqv}{\ensuremath{
\mathchoice{\quad \Longleftrightarrow \quad}{\Leftrightarrow}
                {\Leftrightarrow}{\Leftrightarrow}}}
\newcommand{\imp}{\ensuremath{\Rightarrow}}
\newenvironment{ack}{\medskip{\it Acknowledgement.}}{}

\begin{document}

\authortitle{Anders Bj\"orn and Jana Bj\"orn}
            {Sharp Besov capacity estimates for annuli in metric spaces
with doubling measures}

\author{
Anders Bj\"orn \\
\it\small Department of Mathematics, Link\"oping University, SE-581 83 Link\"oping, Sweden\\
\it \small anders.bjorn@liu.se, ORCID\/\textup{:} 0000-0002-9677-8321
\\
\\
Jana Bj\"orn \\
\it\small Department of Mathematics, Link\"oping University, SE-581 83 Link\"oping, Sweden\\
\it \small jana.bjorn@liu.se, ORCID\/\textup{:} 0000-0002-1238-6751
}

\date{Preliminary version, \today}
\date{}

\maketitle

\noindent{\small
{\bf Abstract}.
We obtain precise estimates, in terms of the measure of balls, 
for the Besov  capacity of annuli and singletons in complete metric spaces.
The spaces are only assumed
to be uniformly perfect with respect to the centre of the annuli
 and equipped with a doubling 
measure.
}

\medskip


\medskip


\section{Introduction}

Capacities are intimately related to function spaces in the sense that 
various properties, such as quasicontinuity and Lebesgue points, 
of functions in such spaces are measured by capacity.
Capacities also reflect metric and measure-theoretic properties 
of the underlying space on which they are defined. 
For example, it is well known that the \p-capacity of a spherical condenser 
in $\R^n$ with
$0<2r\le R$ reflects the dimension of the space as follows,
\begin{equation}  \label{eq-p-cap-Rn}
\cp(B(x,r),B(x,R)) \simeq \begin{cases}
    r^{n-p} &\text{if } 1\le p<n, \\
    (\log(R/r))^{1-p} &\text{if } 1\le p=n, \\
    R^{n-p} &\text{if } p>n. \end{cases}
\end{equation}
Capacities also play an important 
role in fine potential theory and appear 
in the famous Wiener criterion characterizing
boundary regularity for various equations, such as $\Delta_pu=0$
(Maz\cprime ya~\cite{mazya70} and Kilpel\"ainen--Mal\'y~\cite{KilMa-Acta},
with the \p-capacity as in~\eqref{eq-p-cap-Rn}) and
the fractional \p-Laplace equation
$(-\Delta_p)^s u =0$ (Kim--Lee--Lee~\cite{KLL}, using a fractional 
Besov capacity from~\cite{JBWien}).

In this paper we study Besov capacities
on a complete metric space $Y=(Y,d)$ equipped with a doubling measure $\nu$.
Analogously to~\eqref{eq-p-cap-Rn}, we
are primarily interested in estimates for (thick) annuli,
i.e.\ of the capacity for a ball $B(x_0,r)$ within $B(x_0,R)$ where $0 < 2r \le R$.

\medskip

\emph{Throughout the paper we assume that $1  \le p < \infty$.
We also fix a point $x_0$ and let $B_r=B(x_0,r)$ be the open
ball with radius $r$ and centre $x_0$.}

\medskip

The following are our main results.

\begin{thm}\label{thm-metaest} 
Assume that
$Y$ is a complete metric space which
is uniformly perfect at $x_0$ and equipped with a doubling 
measure $\nu$. Let $p>1$ and $0<\theta<1$.
Then for all $0<2r\le R \le \tfrac14\diam Y$, 
\begin{equation}  \label{eq-metaest}
\cpt(\BY_r,\BY_R) 
\simeq \biggl( \int_r^R \biggl( \frac{\rho^{\theta p}}{\nu(\BY_\rho)} \biggr)^{1/(p-1)} 
           \,\frac{d\rho}{\rho} \biggr)^{1-p}
\end{equation}
and
\begin{equation}  \label{eq-metaest-0}
\cpt(\{x_0\},\BY_R) 
\simeq \biggl( \int_0^R \biggl( \frac{\rho^{\theta p}}{\nu(\BY_\rho)} \biggr)^{1/(p-1)} 
           \,\frac{d\rho}{\rho} \biggr)^{1-p},
\end{equation}
with the comparison constants in\/ \textup{``$\simeq$''}
independent of $x_0$, $r$ and $R$.
\end{thm}

Here, $\cpt$ is the Besov condenser capacity defined
for bounded sets $E \Subset \Om$ as
\begin{equation}
\label{eq-Bes-cap-energy}
   \cpt(E,\Om) = \inf_u \int_Y\int_Y \frac{|u(x)-u(y)|^p}{d(x,y)^{\theta p}} 
      \frac{d\nu(y)\, d\nu(x)}{\nu(B(x,d(x,y)))},
\end{equation}
where the infimum is taken over all $u$
such that $0 \le u \le 1$ everywhere,
$u = 1$  in a neighbourhood of $E$ and $\supp u\Subset\Om$, 
see Definition~\ref{def-Besov-condenser-cap}.

The Euclidean spaces and their subsets, 
equipped with the Lebesgue measure or weighted measures $w\,dx$, 
and even singular doubling measures, are included as special cases of our results.
We emphasize that we do not assume any 
Poincar\'e inequalities for upper gradients on $Y$
(as in Gogatishvili--Koskela--Zhou~\cite[Section~4]{GKZ}
and Koskela--Yang--Zhou~\cite{KosYangZhou}).
This makes our results applicable to many disconnected
spaces and spaces carrying few rectifiable curves, including
fractals.

To formulate the next two results we need 
the following exponent sets:
\begin{align*}
  \lQo  &=\biggl\{q>0 : 
        \frac{\nu(B_r)}{\nu(B_R)}  \simle \Bigl(\frac{r}{R}\Bigr)^q 
        \text{ for } 0 < r < R \le 1
        \biggr\}, \\
  \lSo
          &=\{s>0 : 
        \nu(B_r)  \simle r^s 
        \text{ for } 0 < r  \le 1
        \}, \\
  \uSo
          &=\{s>0 : 
        \nu(B_r)  \simge  r^s 
        \text{ for } 0 < r  \le 1
        \}, \\
  \uQo
       &=\biggl\{q>0 : 
       \frac{\nu(B_r)}{\nu(B_R)} 
       \simge  \Bigl(\frac{r}{R}\Bigr)^q 
       \text{ for } 0 < r < R \le 1
       \biggr\}.
\end{align*}
These 
sets were introduced in
Bj\"orn--Bj\"orn--Lehrb\"ack~\cite{BBLeh1}
to capture the local behaviour of the measure at $x_0$.
For example, for the Lebesgue measure in $\R^n$,
\[
\lQo = \lSo = (0,n] \quad \text{and} \quad \uSo = \uQo = [n,\infty).
\]
The subscript $0$ in the above definitions  stands for the fact
that the inequalities are required to hold for small radii.
It is easily verified (see 
\cite[Lemmas~2.4 and~2.5]{BBLeh1})
that the exponent sets can equivalently be defined using
$0<r\le \Theta R\le R_0$ for any fixed $0<\Theta<1$ and $R_0>0$,
even though the implicit comparison constants in ``$\simle$''
and ``$\simge$'' will then depend on $\Theta$ and $R_0$.

All of these sets are intervals. 
The reason for introducing them
as sets is that they may or may not contain their endpoints
\begin{equation} \label{eq-uso}
\quad     \lqo = \sup \lQo,
\quad     \lso = \sup \lSo, 
\quad   \uso = \inf \uSo
\quad \text{and} \quad
\uqo = \inf \uQo.
\end{equation}
Note that $\uqo<\infty$ if 
$\nu$ is doubling, and that
$\lqo>0$ if $Y$ is also uniformly  perfect at $x_0$
(see Heinonen~\cite[Exercise~13.1]{Heinonen}).

When $p>1$ and $\theta p < \lqo$ or $\theta p > \uqo$, 
Theorem~\ref{thm-metaest} provides us with exact estimates for 
the capacity $\cpt(B_r,B_R)$ in terms of $\nu(B_r)$ or $\nu(B_R)$.
When $p=1$, Theorem~\ref{thm-metaest}  cannot be used, 
but we obtain the following similar estimates for $\cpt(\BY_r,\BY_R)$
by using results from Bj\"orn--Bj\"orn--Lehrb\"ack~\cite{BBLeh1},
which cover all $p\ge1$.
The borderline cases $\theta p= \max \lQo$
and $\theta p= \min \uQo$ are considered in
Theorem~\ref{thm-borderline}.
See also Remarks~\ref{rmk-lQ} and~\ref{rmk-dep-const-lqo}.

\begin{thm} \label{thm-lqo}
Assume that $Y$ is a complete metric space which
is uniformly perfect at $x_0$ and equipped with a doubling 
measure $\nu$. Let $0<\theta<1$ and
$0 < R_0 \le \tfrac{1}{4} \diam Y$, with $R_0$ finite.
\begin{enumerate}
\item  \label{it-1}
If $\theta p < \lqo$, then
\begin{equation} \label{eq-it-1}
\cpt(B_r,B_{R})\simeq \frac{\nu(B_r)}{r^{\theta p}}
\quad \text{whenever\/ } 0<2r \le R \le R_0. 
\end{equation}
\item  \label{it-2}
If $\theta p > \uqo$, then
\begin{equation} \label{eq-it-2}
\cpt(B_r,B_{R})\simeq \frac{\nu(B_{R})}{R^{\theta p}}
\quad \text{whenever\/ } 0<2r \le R \le R_0.
\end{equation}
\end{enumerate}
\bigskip
%

\noindent
In both cases, the comparison constants in\/ \textup{``$\simeq$''}
depend on $R_0$.

Moreover, the lower bound in \eqref{eq-it-1} implies $p \in \lQo$,
while the lower bound in \eqref{eq-it-2} implies $p \in \uQo$.
If $p>1$ then \eqref{eq-it-2} holds if and only if
$\theta p > \uqo$.
\end{thm}

In Ahlfors regular spaces, estimates \eqref{eq-it-1}
and \eqref{eq-it-2} were given in Lehrb\"ack--Shan\-mu\-ga\-lin\-gam~\cite{LS22},
and used to show that Besov-norm-preserving homeomorphisms 
between such spaces are  quasisymmetric.

In many situations it is important whether singletons have
zero or positive capacity. 
In the following result, we characterize these cases in terms of
the exponent sets $\uSo$ and $\lSo$.

\begin{thm} \label{thm-S}
Assume that $Y$ is a complete metric space which
is uniformly perfect at $x_0$ and equipped with a doubling 
measure $\nu$. Let $0<\theta<1$.
\begin{enumerate}
\item  \label{it-1-s}
  If $\theta p > \uso$,
  then
\[
\Cpt(\{x_0\})>0
\quad \text{and} \quad
\cpt(\{x_0\},B_R)>0 \text{ for every }0 < R <\tfrac12 \diam Y,
\]
where the capacity $\Cpt$ is defined by means of
  the Besov norm as in Definition~\ref{eq-def-norm-cap}.
\item  \label{it-2-s}
  If $\theta p \notin \uSo$
\textup{(}in particular if $\theta p < \uso$\/\textup{)},
  or if $p>1$ and $\theta p \in \lSo$, then
\[
\Cpt(\{x_0\})=0
\quad \text{and} \quad
\cpt(\{x_0\},B_R)=0 \text{ for every }R>0.
\]
\end{enumerate}
\end{thm}  
\bigskip

In Anttila~\cite{anttila}, the numbers $\us_0$ and $\ls_0$ are called
the  \emph{upper and lower local dimensions} of $\mu$ at $x_0$,
while $\uq$ in Remark~\ref{rmk-lQ} is called
the \emph{pointwise Assouad dimension} of $\mu$ at $x_0$.
(See~\cite[Lemma~2.4]{BBLeh1} for why
the definitions of $\us_0$ and $\ls_0$ in~\cite{anttila}
are equivalent to those in~\eqref{eq-uso}.)
In \cite{BBLehIntGreen}, $\us_0$ played a decisive role in determining
the sharp integrability properties for \p-harmonic Green functions and
their gradients.

On $\R^n$, the spaces defined by means of the energy integral 
 in~\eqref{eq-Bes-cap-energy}
are often called fractional Sobolev spaces and are the traces of Sobolev spaces 
on sufficiently nice domains
(Jonsson--Wallin~\cite{JW84}).
As such, they are suitable as boundary values for various Dirichlet problems and 
appear in boundary regularity results for elliptic differential
equations (Kristensen--Mingione~\cite{KriMin10}).

They also play an important role 
in nonlocal problems, such as the fractional \p-Laplace equation
$(-\Delta_p)^s u =0$.
These problems have attracted a lot of attention in the past two decades, see e.g.\ 
Kim--Lee--Lee~\cite{KLL},
Korvenpää--Kuusi--Lindgren~\cite{KorKuuLin} and 
Lindgren--Lindqvist~\cite{LinLinEigenv},
to name just a few.

Recently, similar problems and the associated Besov spaces have been
studied for metric measure spaces in e.g.\ 
Capogna--Kline--Korte--Shan\-mu\-ga\-lin\-gam--Snipes~\cite{capognaKKSS},
Eriksson-Bique--Giovannardi--Korte--Shan\-mu\-ga\-lin\-gam--Speight~\cite{EBGKSS},
Go\-ga\-tish\-vi\-li--Koskela--Shan\-mu\-ga\-lin\-gam~\cite{GKS},
Go\-ga\-tish\-vi\-li--Koskela--Zhou~\cite{GKZ}
and Koskela--Yang--Zhou~\cite{KosYangZhou}.
The role of Besov spaces as traces of Sobolev type spaces was in the 
metric setting studied in
Bourdon--Pajot~\cite{BP},
Bj\"orn--Bj\"orn--Gill--Shan\-mu\-ga\-lin\-gam~\cite{BBGS}
and
Bj\"orn--Bj\"orn--Shan\-mu\-ga\-lin\-gam~\cite{BBShyptrace}, 
and will be one of our main tools.

Our approach to the above estimates is based on extensions of Besov
functions from $Y$ to hyperbolic fillings of $Y$, together with 
estimates from Bj\"orn--Bj\"orn--Lehrb\"ack~\cite{BBLeh1}, \cite{BBLehIntGreen}
for \p-capacities associated with Sobolev spaces.
More precisely, we use the comparison between Besov seminorms of functions on $Y$
and the Dirichlet energy of their extensions to a uniformized
hyperbolic filling of $Y$, obtained in~\cite{BBShyptrace}.
These constructions and comparisons are done in Section~\ref{sect-hypfill}.

However, since the results in~\cite{BBShyptrace} only cover bounded 
spaces, special care has to be taken for unbounded $Y$.
This is done in Section~\ref{sect-unbdd} by replacing $Y$ with a 
suitably chosen bounded subset, so that the restriction of $\nu$ is still doubling.
Even when $Y$ is bounded, it is only biLipschitz equivalent to the 
boundary  of the uniformized hyperbolic filling of $Y$, 
which would in turn put serious restrictions on the allowed radii
$r$ and $R$ in our estimates. 
In Section~\ref{sect-ext-Y} we therefore show how to replace $Y$
by a carefully constructed enlarged space so that the involved
capacities are comparable and all radii $\le\tfrac14\diam Y$ can be treated.

Along the way, in Sections~\ref{sect-Besov} and~\ref{sect-doubl-cap},
we prove various fundamental properties of Besov capacities
in metric spaces, both for doubling and nondoubling measures,
including in some cases also $\theta\ge1$.
Finally, in Sections~\ref{sect-sharp-est} and~\ref{sect-borderline},
we prove Theorems~\ref{thm-metaest}--\ref{thm-S}.

As mentioned above, we use hyperbolic fillings to obtain our main results.
  It would be interesting to find more direct proofs.
On the other hand, our technique shows that there is a direct
  correspondence between these results and 
the corresponding
results for Sobolev spaces
in~\cite{BBLeh1} and~\cite{BBLehIntGreen}.

\begin{ack}
A.~B. and J.~B. were supported by the Swedish Research Council,
  grants 2020-04011 resp.\ 2018-04106.
Part of this research was done 
when the authors
visited Institut Mittag-Leffler in the autumn of 2022 during the programme
\emph{Geometric Aspects of Nonlinear Partial Differential Equations}.
We thank the institute for their hospitality
and support.
\end{ack}

\section{Preliminaries}
\label{sect-prelim}

\emph{In this section we assume that 
$X=(X,d)$ is a metric space equipped
with a Borel measure $\mu$ such that $0 < \mu(B)<\infty$ for every ball $B \subset X$.}
To avoid pathological situations we also assume
that all metric spaces, considered
in this paper,
contain at least two points.

\medskip

As is often customary we extend $\mu$, and other measures, as
  outer measures defined on all sets.
  This plays a role at least  in
  Proposition~\ref{prop-Cpt-properties}\,\ref{Cp-mu<=Cp}.

A metric space is \emph{proper} if all closed bounded sets are compact.
We denote balls in $X$ by 
\[
B(x,r)=\{y \in X: d(y,x) <r\} \quad \text{and let 
$\la B(x,r)=B(x,\la r)$.}
\]
All balls in this paper are open.
In metric spaces it can happen that
balls with different centres and/or radii denote the same set. 
We will however use the convention that a ball comes with
a predetermined centre and radius.

The space $X$ is \emph{uniformly perfect at $x$} if there
is a constant $\ka>1$ such that
\begin{equation} \label{eq-unif-perfect}
  B(x,\ka r) \setm B(x,r) \ne \emptyset
\quad \text{whenever } B(x,\ka r) \ne X.
\end{equation}  
In fact, it then follows that \eqref{eq-unif-perfect} holds
whenever $B(x,r) \ne X$,
since if
$B(x,\ka r) = X$ then $B(x,\ka r) \setm B(x,r) = X \setm B(x,r) \ne \emptyset$.

The space $X$ is \emph{uniformly perfect} if it is uniformly
perfect at every $x$ with the same constant $\ka$.
This definition coincides with the one 
in Heinonen~\cite[Section~11.1]{Heinonen}, see therein for more on the
history of this assumption.
We do not know if pointwise uniform perfectness has been used before.
Note that $X$ is uniformly perfect with any $\ka>1$ if $X$ is connected.

The measure $\mu$  is \emph{doubling}
if  there is a \emph{doubling constant} $C_\mu> 1$ such that 
\[ 
  0 < \mu(2B)\le C_\mu \mu(B) < \infty
\quad \text{for all balls $B$}.
\]
Similarly, $\mu$ is \emph{reverse-doubling at $x$}, if there
are constants $C,\kard>1$ such that 
\begin{equation}\label{eq:rev-doubling-x0}
  \mu(B(x,\kard r))\ge C \mu(B(x,r))
\quad \text{for all $0<r < \diam X/2\kard$}.
\end{equation}
By continuity of the measure,
the estimate \eqref{eq:rev-doubling-x0}
holds also if $r = \diam X/2\kard< \infty$, as required
in Bj\"orn--Bj\"orn--Lehrb\"ack~\cite{BBLeh1}.
If $\mu$ is doubling, it is easy to see that 
$X$ is uniformly perfect at $x$
if and only if $\mu$ is reverse-doubling at $x$.
(For necessity we can choose any $\kard > \ka$,
and for sufficiency any $\ka > 2\kard$.)
If $\mu$ is doubling and $X$ is connected,
then $\mu$ is reverse-doubling at every $x$ with any $\kard>1$.

Throughout the paper, we write $a \simle b$ if there is an implicit
constant $C>0$ such that $a \le Cb$, and analogously $a \simge b$ if $b \simle a$,
and $a \simeq b$ if $a \simle b \simle a$.
The implicit comparison constants are allowed to depend on the standard parameters.
We will carefully explain the dependence in each case.
See Remarks~\ref{rmk-thm-1.1} and~\ref{rmk-dep-const-lqo}
for the dependence in Theorems~\ref{thm-metaest} and~\ref{thm-lqo}.

Sometimes, when dealing with several different spaces simultaneously,
we will write $B^X$,  $d_X$, $\uQo^X$, $\uqo^X$ etc.\ to indicate that these notions are taken with respect
to the metric space $X$.

\section{Besov spaces and capacities}
\label{sect-Besov}

\emph{In this section we assume that $Y=(Y,d)$ is a
proper metric space equipped
with a Borel measure~$\nu$ such that $0 < \nu(B) <\infty$ for every ball
  $B \subset Y$.
  We also assume that $\theta>0$, and emphasize that in this section
  $\theta \ge 1$ is allowed.}

\medskip

For a measurable $u: Y \to [-\infty,\infty]$ (which is
finite $\nu$-a.e.) we define the \emph{Besov seminorm} by
\begin{equation*} 
  [u]_{\theta,p}=
[u]_{\theta,p,Y}=
\biggl(  \int_Y\int_Y \frac{|u(x)-u(y)|^p}{d(x,y)^{\theta p}} 
      \frac{d\nu(y)\, d\nu(x)}{\nu(B(x,d(x,y)))} \biggr)^{1/p}. 
\end{equation*}
Here and elsewhere, the  integrand 
should be interpreted as zero  when $y=x$.

The Besov space $B^\theta_p(Y)$ consists of the functions $u$ such that
the \emph{Besov norm}
\begin{equation}   \label{eq-def-Besov-spc}
\|u\|_{B^\theta_p(Y)}^p:=[u]_{\theta,p}^p+\|u\|_{L^p(Y)}^p < \infty.
\end{equation}
This space is a Banach space,
see Remark~9.8 in Bj\"orn--Bj\"orn--Shan\-mu\-ga\-lin\-gam~\cite{BBShyptrace}.
(The norm \eqref{eq-def-Besov-spc} is equivalent to the one in \cite{BBShyptrace},
but
the norm-capacity $\Cpt$ below exactly coincides with the one in \cite{BBShyptrace}.)

We restrict our attention to Besov spaces with two indices
(i.e.\ ``$q=p$''). 
Such Besov spaces are often called fractional Sobolev spaces or
  Sobolev--{Slobodetski\u\i} spaces, although Besov spaces seem
  to be the most common name in the metric space literature.

Assuming that $\nu$ is doubling,
equivalent definitions, using equivalent seminorms, are given in
Gogatishvili--Koskela--Shan\-mu\-ga\-lin\-gam~\cite[Theorem~5.2 and (5.1)]{GKS}.
When $\nu$ is also reverse-doubling (or equivalently, uniformly perfect),
further equivalent definitions can be found in
Gogatishvili--Koskela--Zhou~\cite[Theorem~4.1 and Proposition~4.1]{GKZ}, 
for example that
the Besov space $B^\theta_p(Y)$ considered here
coincides with
the corresponding Haj\l asz--Besov space.
By \cite[Lemmas~6.1 and~6.2]{GKS}, it is also related to fractional
Haj\l asz spaces, considered already in Yang~\cite{Yang2003}.
See these papers for the precise definitions and
earlier references to the theory on
$\R^n$ and on Ahlfors regular metric spaces.

We are interested in two types of Besov capacities, the norm-capacity and
the condenser capacity.

\begin{deff}
The \emph{Besov norm-capacity} of $E \subset Y$ is
\begin{equation}  \label{eq-def-norm-cap}
\Cpt(E) = \CptY(E)
=\inf_u {\|u\|_{B^\theta_p(Y)}^p},
\end{equation}
where the infimum is taken over all $u\in B_p^\theta(Y)$ such that 
$u = 1$  in a neighbourhood of $E$ and $0 \le u \le 1$ everywhere.
Such $u$ are called \emph{admissible} for $\Cpt(E)$.
\end{deff}

By truncation it follows that one can equivalently take
the infimum over all $u$ such that $u \ge 1$ in
a neighbourhood of $E$.
As usual, when requiring that $u \ge 1$ or $0 \le u \le 1$ everywhere
we mean that there is a representative of $u$ 
satisfying these requirements.
By $E \Subset \Om$ we mean that $\itoverline{E}$
is a compact subset of $\Om$.

\begin{deff}  \label{def-Besov-condenser-cap}
Let $\Om \subset Y$ be a bounded
open set and $E \Subset \Om$.
The \emph{Besov condenser capacity} is given by
\[
\cpt(E,\Om) = \cptY(E,\Om) = \inf_u {[u]_{\theta,p}^p},
\]
where the infimum is taken over all $u$
such that $0 \le u \le 1$ everywhere,
$u = 1$  in a neighbourhood of $E$
and $\supp u\Subset\Om$.
Such $u$ are called \emph{admissible} for $\cpt(E)$.
\end{deff}

The corresponding capacities for Sobolev spaces are called Sobolev resp.\
variational capacity in \cite{BBbook}.
Condenser capacities are also often called ``relative''.

There do not seem to be very many papers on Besov capacities in metric spaces.
Nuutinen~\cite{Nuutinen}
and Heikkinen--Koskela--Tuominen~\cite{HKT2017}
extensively studied the norm-capacity, defined using
the Haj\l asz--Besov norm, under
the assumption that $\nu$ is doubling.
In~\cite{HKT2017},
they also considered the corresponding
Triebel--Lizorkin norm-capacity, which was later studied
by Karak~\cite{karak20}.
The Besov norm-capacity  $\Cpt$ was used by
Bj\"orn--Bj\"orn--Shan\-mu\-ga\-lin\-gam~\cite{BBShyptrace}.
The Besov condenser capacity $\cpt$ was studied
in the Ahlfors $Q$-regular case 
by Bourdon~\cite{bourdon07} ($p >Q$ and $\theta=1/p$),
Costea~\cite{costea09publmat} ($p >Q$) and
Lehrb\"ack--Shan\-mu\-ga\-lin\-gam~\cite{LS22}.

Our main estimates remain the same (up to changes in implicit constants)
when the seminorm is replaced by an equivalent seminorm.
However, some of the basic properties, such as subadditivity,
 are not directly transferable, although
the proofs often are, so we include them here.

\begin{prop} \label{prop-Cpt-properties}
Let $E,E_1,E_2,\ldots \subset Y$.
Then the following properties hold\/\textup{:}
\begin{enumerate}
\renewcommand{\theenumi}{\textup{(\roman{enumi})}}%
\item \label{Cp-subset-sum}
  if $E_1 \subset E_2$, then $\Cpt(E_1) \le \Cpt(E_2)$,
\item \label{Cp-mu<=Cp}
  $\nu(E) \le \Cpt(E)$,
\item \label{Cp-Choq-K-sum}
if $K_1 \supset K_2 \supset \cdots$ are compact subsets of $X$, then
\[
      \Cpt\biggl(\bigcap_{i=1}^\infty K_i\biggr) 
          = \lim_{i \to \infty} \Cpt(K_i),
\]
\item \label{Cp-subadd-sum}
  $\Cpt$ is countably subadditive, i.e.\
  if $E=\bigcup_{i=1}^\infty E_i$ then
\[
      \Cpt(E) 
          \le \sum_{i=1}^\infty \Cpt(E_i).
\]
\end{enumerate}  
\end{prop}

The monotonicity~\ref{Cp-subset-sum} is trivial,
while \ref{Cp-mu<=Cp} follows directly from the definition.
The property~\ref{Cp-Choq-K-sum}
follows from the fact that $\Cpt$ is an outer capacity (by definition), i.e.
\[
         \Cpt(E) = \inf_{\substack{ G \supset E \\ G \text{ open}}} \Cpt(G),
\]
and elementary properties of compact sets, see 
Nuutinen~\cite[Section~3]{Nuutinen}.
As for~\ref{Cp-subadd-sum}, Nuutinen~\cite{Nuutinen} only obtains  quasi-subadditivity
since he works in a more general setting
in which the countable subadditivity does not always hold.
We therefore provide a proof.

\begin{proof}[Proof of \ref{Cp-subadd-sum}]
We may assume that the right-hand side is finite.
Let $\eps>0$.
For each $i=1,2,\ldots$\,, choose
$u_i$ admissible for $\Cpt(E_i)$ with
\[
  [u_i]_{\theta,p}^p+\|u_i\|_{L^p(Y)}^p < \Cpt(E_i)+ \frac{\eps}{2^i}.
\]
Let $u=\sup_i u_i$.
Then $u = 1$ in a neighbourhood of $\bigcup_{i=1}^\infty E_i$.
Moreover, for $x,y \in Y$, 
\[
  |u(x)-u(y)|^p
  \le \sup_i |u_i(x)-u_i(y)|^p
  \le \sum_{i=1}^\infty |u_i(x)-u_i(y)|^p
\]
and similarly,
$|u(x)|^p =   \sup_i |u_i(x)|^p 
  \le \sum_{i=1}^\infty |u_i(x)|^p$.
Hence
\begin{align*}
      \Cpt(E) 
      & \le ([u]_{\theta,p}^p+\|u\|_{L^p(Y)}^p)
       \le \sum_{i=1}^\infty ([u_i]_{\theta,p}^p+\|u_i\|_{L^p(Y)}^p) \\
      & 
           < \sum_{i=1}^\infty \Bigl( \Cpt(E_i)+ \frac{\eps}{2^i}\Bigr) 
           =   \sum_{i=1}^\infty \Cpt(E_i) + \eps.
\end{align*}
Letting $\eps \to 0$ completes the proof.
\end{proof}  

\begin{prop} \label{prop-cp}
Let $\Om \subset \Om' \subset Y$ be  open and bounded,
$E, E_1, E_2,\ldots \Subset \Om$. 
Then the following properties hold\/\textup{:}
\begin{enumerate}
\renewcommand{\theenumi}{\textup{(\roman{enumi})}}%
\item \label{cp-subset-sum}
  if $E_1 \subset E_2 \subset \Om$,
 then\/ $\cpt(E_1,\Om') \le \cpt(E_2,\Om)$,
\item \label{cp-Choq-K-sum}
if $K_1 \supset K_2 \supset \cdots $ are compact subsets 
of\/ $\Om$, then
\[
      \cpt\biggl(\bigcap_{i=1}^\infty K_i,\Om\biggr) 
          = \lim_{i \to \infty} \cpt(K_i,\Om),
\]
\item \label{cp-subadd-sum}
$\cpt$ is countably subadditive, i.e.\ if $E=\bigcup_{i=1}^\infty E_i$ then
\[
      \cpt(E)
          \le \sum_{i=1}^\infty \cpt(E_i,\Om).
\]
\end{enumerate}
\end{prop}

Again, \ref{cp-subset-sum} is trivial,
while \ref{cp-Choq-K-sum}
follows from elementary
properties of compact sets  since $\cpt$ is
an outer capacity (by definition).
The proof of \ref{cp-subadd-sum} is similar to the proof of
Proposition~\ref{prop-Cpt-properties}\,\ref{Cp-subadd-sum}.

In the Ahlfors $Q$-regular case with $p >Q > 1$,
these facts were stated in Costea~\cite{costea09publmat}
with a comment that the proof is essentially the same as 
in Costea~\cite[Theorem~3.1]{costea07}.
His proof of \ref{cp-subadd-sum}
uses reflexivity.
Our proof is considerably shorter and also covers 
the case $1\le p\le Q$
as well as the non-Ahlfors regular case.

\section{Capacity estimates when
  \texorpdfstring{$\nu$}{nu} is doubling}

\label{sect-doubl-cap}

\emph{In this section we assume that
$Y$ is a complete metric space equipped
with a doubling measure~$\nu$ and that $0<\theta<1$.}

\medskip

Note that $Y$ is proper,
  see Bj\"orn--Bj\"orn~\cite[Proposition~3.1]{BBbook}.
The comparison constants in this section are independent of the choice of $x_0$,
they depend only on $\theta$, $p$  
  and $C_\nu$,
  unless said otherwise.

Our next aim is to deduce the following result,
  which will be important later on.
  Note that $B_{2R} \ne Y$ whenever $R < \tfrac14 \diam Y$.

\begin{prop} \label{prop-cpY-Theta}
Assume that $Y$ is uniformly perfect at $x_0$ with constant $\ka$.
  Fix   $0 < \Theta < 1$.
  If  $0 < \Theta R \le  2r  \le R$ and
 $B_{2R}  \ne Y$, then
  \[
  \cpt(B_r,B_R)
      \simeq \frac{\nu(B_r)}{r^{\theta p}}
      \simeq \frac{\nu(B_R)}{R^{\theta p}},
   \]
   with comparison constants also depending on $\ka$ and $\Theta$.
\end{prop}  

We split the proof of Proposition~\ref{prop-cpY-Theta} into
two parts.
We begin with the lower bound, which holds also when $\theta \ge 1$.

\begin{prop} \label{prop-cpt-ge}
  Assume that $Y$ is uniformly perfect at $x_0$ with constant $\ka$, and that
$0 < 2r  \le R $ with $B_{2R} \ne Y$.
  Then
  \[
   \cpt(B_r,B_R) \simge \frac{\nu(B_r)}{R^{\theta p}},
   \]
with comparison constant also depending on $\ka$.
\end{prop}

\begin{proof}
Let $u$ be admissible for $\cpt(B_r,B_R)$.
As $B_{2R}\ne Y$ it follows from the uniform perfectness that there exists
$z\in B_{2\ka R} \setm B_{2R}$.
Since $B(z,R) \cap B_R = \emptyset$ and
$d(x,y)\le (2\ka+2)R$ for all $x\in B(z,R)$ and $y\in B_r$, we get that
\[ 
  [u]_{\theta,p}^p
  \ge
   \int_{B(z,R)}\int_{B_r}   \frac{1}{((2\ka+2)R)^{\theta p}}  
      \frac{d\nu(y)\, d\nu(x)}{\nu(B(x,(2\ka+2)R))} 
  \simge
  \frac{\nu(B_r)}{R^{\theta p}}.
\] 
Taking infimum over all $u$ admissible for $\cpt(B_r,B_R)$ concludes the proof.
\end{proof}

To prove the upper bound in Proposition~\ref{prop-cpY-Theta}
  we will use the following  simple lemma,
  which will also be used when proving   Lemma~\ref{lem-cp-Cp}.

\begin{lem} \label{lem-Lip}
  Assume   that $0\le\eta\le1$ is an $M$-Lipschitz function on $Y$.
If $x \in Y$, then
\[
  I(x):= \int_Y  \frac{|\eta(x)-\eta(y)|^p}{d(x,y)^{\theta p}} 
   \frac{d\nu(y)}{\nu(B(x,d(x,y)))}
   \simle M^{\theta p}.
\]
\end{lem}  

\begin{proof}
Let $B^j=B(x,2^j/M)$.
Since $\nu(B(x,d(x,y))) \simeq \nu(B^{j})$
for $y\in B^{j} \setm B^{j-1}$ 
  and $0<\theta<1$, we see that
\begin{align*}
   I(x) &\simeq   \sum_{j=-\infty}^\infty
  \int_{B^j \setm B^{j-1}}  \frac{|\eta(x)-\eta(y)|^p}{d(x,y)^{\theta p}} 
  \frac{d\nu(y)}   {\nu(B^{j})}  \\   
  &  \simle \sum_{j=-\infty}^0  \int_{B^{j} \setm B^{j-1}}
  \frac{M^p d(x,y)^{(1-\theta) p}\, d\nu(y)} {\nu(B^{j})} 
+ \sum_{j=1}^\infty  \int_{B^{j} \setm B^{j-1}}
  \frac{d(x,y)^{-\theta p} \,d\nu(y)} {\nu(B^{j})}  \\
 & \simle \sum_{j=-\infty}^0 M^{\theta p} 2^{j(1-\theta) p} 
+  \sum_{j=1}^\infty M^{\theta p} 2^{-j\theta p}
\simeq M^{\theta p}.
\qedhere
\end{align*}
\end{proof}

This now leads to the following estimate.

\begin{prop} \label{prop-cpt-le-new}
  Assume that   $0 < 2r  \le R$.
  Then
  \[
  \cpt(B_r,B_R) \simle
  \min\biggl\{\frac{\nu(B_r)}{r^{\theta p}},\frac{\nu(B_R)}{R^{\theta p}}\biggr\}.
   \]
\end{prop}

\begin{proof}
Let  $u$ be a $3/R$-Lipschitz function admissible
for $\cpt(B_{R/2},B_R)$.
The doubling property and symmetry in $x$ and $y$ imply that
\begin{equation*}   
[u]_{\theta,p}^p \simeq \int_{B_R} \int_{Y}
\frac{|u(x)-u(y)|^p}{d(x,y)^{\theta p}} 
\frac{d\nu(y)}{\nu(B(x,d(x,y)))} \,d\nu(x).
\end{equation*}
Integrating the estimate from Lemma~\ref{lem-Lip} over $x\in B_R$
gives
\[
\cpt(B_r,B_R) \le \cpt(B_{R/2},B_R)
\simle   \frac{\nu(B_R)}{R^{\theta p}}.
\]
Applying the last estimate
 with $R$ replaced by $2r$ gives
\[
\cpt(B_r,B_R) \le \cpt(B_{r},B_{2r})
\simle   \frac{\nu(B_r)}{r^{\theta p}}. \qedhere
\]
\end{proof}

\begin{proof}[Proof of Proposition~\ref{prop-cpY-Theta}]
  This follows directly from
  Propositions~\ref{prop-cpt-ge} and~\ref{prop-cpt-le-new},
  together with the doubling property.
\end{proof}

\begin{lem} \label{lem-cp-Cp}
Assume that $\Om \subset Y$ is open and bounded,
and let $E  \Subset \Om$.
If $\Cpt(E)=0$, then $\cpt(E,\Om)=0$.
\end{lem}

\begin{proof}
Since $E \Subset \Om$ there is a Lipschitz function $0 \le \eta \le 1$
such that $\eta =1$ in a neighbourhood of $E$ and $\eta=0$
in a neighbourhood of $Y \setm \Om$.
Let $M$ be the Lipschitz constant of $\eta$ and let 
$\eps>0$.
As $\Cpt(E)=0$ there is a function  $u$ admissible
for $\Cpt(E)$ with 
$\|u\|_{B^\theta_p(Y)}^p < \eps$.
Let $v=u\eta$.
Then 
\begin{align*}
|v(x)-v(y)|
& = |u(x)\eta(x)-u(x)\eta(y)+u(x)\eta(y)- u(y)\eta(y)| \\
& \le u(x) |\eta(x)-\eta(y)| + |u(x)- u(y)|.
\end{align*}
Hence, by Lemma~\ref{lem-Lip},
\begin{align*}
[v]_{\theta,p}^p
& \le  2^p \int_Y u(x)^p
\int_Y  \frac{|\eta(x)-\eta(y)|^p}{d(x,y)^{\theta p}} 
\frac{d\nu(y)}{\nu(B(x,d(x,y)))}
+ 2^p [u]_{\theta,p}^p 
\\
   &   \simle  M^{\theta p} \|u\|_{L^p(Y)}^p + [u]_{\theta,p}^p 
  \le (M^{\theta p}+1)  \|u\|_{B^\theta_p(Y)}^p
   < (M^{\theta p}+1)  \eps.
\end{align*}

As $v =1$ in a neighbourhood of $E$ and $v=0$
in a neighbourhood of $Y \setm \Om$, we see that
\[
   \cpt(E,\Om) \le [v]_{\theta,p}^p
   \simle 
   (M^{\theta p}+1)  \eps.
\]
Letting $\eps \to 0$ completes the proof.
\end{proof}

Note that the converse of Lemma~\ref{lem-cp-Cp}
does not hold in general;
consider e.g.\ a compact $Y$ in which case $\cpt(Y,Y)=0$
(as $u \equiv 1$ is admissible) while $\Cpt(Y) \ge \nu(Y) >0$.
Nevertheless, we will prove the following characterization.

\begin{prop} \label{prop-cp-Cp-iff}
  Assume that $\Om \subset Y$ is open and bounded and such that $\nu(Y \setm\Om)>0$.
Let $E  \Subset \Om$.
Then $\Cpt(E)=0$ if and only if $\cpt(E,\Om)=0$.
\end{prop}  

The following simple observation will serve as a Poincar\'e type inequality.
We will use it to prove Proposition~\ref{prop-cp-Cp-iff} as well as
Lemma~\ref{lem-6-16} below.

\begin{lem} \label{lem-quasi-PI}
If $u=0$ outside  a bounded measurable set $\Om$ and 
 $K \subset Y \setm\Om$ is a bounded measurable set with $\nu(K)>0$, then for 
every $z\in K$,
  \[
 \int_Y |u|^p \, d\nu  \le R^{\theta p} \frac{\nu(B(z,R))}{\nu(K)}
         [u]_{\theta,p}^p,
  \]
where 
\[
R=\diam K + \sup \{d(x,y) : x \in K \text{ and } y \in \Om\}.
\]
\end{lem}

\begin{proof}
Since $u=0$
outside $\Om$, and in particular in $K$,
and $B(x,d(x,y))\subset B(z,R)$ for all
$x\in K$ and $y\in\Om$, we see that
\begin{align*}
  \int_Y |u|^p \, d\nu
  & =  \frac{1}{\nu(K)}
     \int_{K}\int_{\Om} |u(x)-u(y)|^p \, d\nu(y)\, d\nu(x) \\
      & \le R^{\theta p} \frac{\nu(B(z,R))}{\nu(K)} 
      \int_K \int_\Om  \frac{|u(x)-u(y)|^p}{d(x,y)^{\theta p}} 
\frac{d\nu(y)}{\nu(B(x,d(x,y)))}.\qedhere
\end{align*}
\end{proof}

\begin{proof}[Proof of Proposition~\ref{prop-cp-Cp-iff}]
One implication was shown in Lemma~\ref{lem-cp-Cp},
so assume that $\cpt(E,\Om)=0$.
By assumption,  there is a compact set $K \subset Y \setm \Om$ with
$\nu(K)>0$.
Let $z \in K$ and  $\eps >0$.
Then there is $u$ admissible for $\cpt(E,\Om)$ with
$[u]_{\theta,p}^p < \eps$.
Since $u$ is admissible also for $\Cpt(E)$,
Lemma~\ref{lem-quasi-PI} implies that
\[
\Cpt(E) \le \biggl( 1+R^{\theta p} \frac{\nu(B(z,R))}{\nu(K)} \biggr) \eps,
\]
and letting $\eps \to 0$ gives $\Cpt(E)=0$.
\end{proof}

As an immediate consequence of Proposition~\ref{prop-cp-Cp-iff}
and monotonicity, we obtain the following characterization.

\begin{cor} \label{cor-singleton}
The following are equivalent\/\textup{:}
\begin{enumerate}
\item \label{a1}
  $\Cpt(\{x_0\})=0$,
\item \label{a2}
  $\cpt(\{x_0\},B_r)=0$ for every $r>0$,
\item \label{a3}
  $\cpt(\{x_0\},B_r)=0$ for some $0 < r < \tfrac12 \diam Y$.
\end{enumerate}
\end{cor}  

If $Y=[-1,1]$ (with Lebesgue measure), $x_0=0$ and $r > 1 = \tfrac12 \diam Y$,
then $u \equiv 1$ is admissible for $\cpt(\{x_0\},B_r)$ and thus
$\cpt(\{x_0\},B_r)=0$.
On the other hand $\Cpt(\{x_0\})>0$ if $\theta p > 1$, by
Theorem~\ref{thm-S}.
This shows that the range in \ref{a3} is sharp.

When $\Om$ is a ball, the following result gives 
more precise 
information than Lemma~\ref{lem-cp-Cp}.

\begin{lem} \label{lem-6-16}
Assume that $E\subset B_r$. Then
\[
\cpt(E,B_{2r}) \simle (1+r^{-\theta p}) \Cpt(E).
\]
If, moreover, $Y$ is uniformly perfect at $x_0$ with constant $\ka$
and $Y\setm B_{3\ka r}\ne\emptyset$, then 
\[
\Cpt(E) \simle (1+r^{\theta p}) \cpt(E,B_{2r}),
\]
with comparison constant also depending on $\ka$.
\end{lem}

\begin{proof}
  Let $u$ be admissible for $\Cpt(E)$ and let $0 \le \eta \le 1$ be
  a $(2/r)$-Lipschitz function such that $\eta =1$ in a neighbourhood of $B_r$
  and $\eta=0$
in a neighbourhood of  $Y\setm B_{2r}$.
Let $v=u\eta$. Then $v$ is admissible for $\cpt(E,B_{2r})$ and 
as in the proof of  Lemma~\ref{lem-cp-Cp},
\[
|v(x)-v(y)| \le u(x) |\eta (x)-\eta(y)| + |u(x)-u(y)|.
\]
Hence by symmetry and Lemma~\ref{lem-Lip},
\begin{align*}
\cpt(E,B_{2r}) \le [v]_{\theta,p}^p &\simle 
\int_{B_{2r}} \int_Y \frac{|v(x)-v(y)|^p}{d(x,y)^{\theta p}} 
      \frac{d\nu(y)\, d\nu(x)}{\nu(B(x,d(x,y)))}  \\
&\simle   r^{-\theta p}  \int_{B_{2r}} |u|^{p}\, d\nu +  [u]_{\theta,p}^p. 
\end{align*}
Taking infimum over all $u$ admissible for $\Cpt(E)$  proves the first
inequality in the statement of the lemma.

For the second inequality, note that every $u$ admissible for $\cpt(E,B_{2r})$
is admissible also for $\Cpt(E)$.
Next, use the uniform perfectness  at $x_0$ to find
$z \in B_{3\ka r} \setm B_{3r}$.
Lemma~\ref{lem-quasi-PI} with $\Om=B_{2r}$, $K=B(z,r)$ 
and $R= (3\ka+3) r$, together with $\nu(B(z,R))\simle \nu(B(z,r))$, then implies that 
\[
\Cpt(E) \le \int_Y |u|^p \, d\nu +  [u]_{\theta,p}^p \simle (1+r^{\theta p} ) [u]_{\theta,p}^p.
\]
Taking infimum over all $u$ admissible for $\cpt(E,B_{2r})$
concludes the proof.
\end{proof}

We conclude this section by comparing capacities
with respect to different underlying spaces.
Since the seminorm $[u]_{\theta,p}$ is nonlocal, the sets 
where $u$ vanishes cannot be ignored.

\begin{lem} \label{lem-cp-restrict-new}
Let $X\subset Y$ be compact,
$\Om \subset X$ be open and   $E \Subset \Om$.
Assume that
  \begin{equation} \label{eq-restrict}
  \nu(B(x,r) \cap X) \simeq \nu(B(x,r))
  \quad \text{for all $x \in X$ and $0 < r < 2 \diam X$},
  \end{equation}
and that for a.e.\ $x\in\Om$,
 \begin{equation}  \label{eq-ass-Y-X-Om}
   \int_{Y \setm X}  I(x,y) \,d\nu(y)
    \simle   \int_{X \setm \Om} I(x,y) \,d\nu(y),
 \end{equation}
 where
 \[
   I(x,y)=\frac{1}{d(x,y)^{\theta p}\,\nu(B(x,d(x,y)))}.
\]

Then 
\[
\cptX(E,\Om) \simeq \cptY(E,\Om)
\]
 with comparison constants also depending on the 
implicit comparison constants in \eqref{eq-restrict} and \eqref{eq-ass-Y-X-Om}.
\end{lem}

\begin{proof}
Note that $u$ is admissible for $\cptX(E,\Om)$ if and only if
its zero extension to $Y\setm X$ is admissible
for $\cptY(E,\Om)$.
Hence it is enough to show that $[u]_{\theta,p,X} \simeq [u]_{\theta,p,Y}$
for any $u$ admissible for $\cptY(E,\Om)$.
Consider such a function $u$.

By \eqref{eq-restrict},
$[u]_{\theta,p,X} \simle [u]_{\theta,p,Y}$.
Conversely, the doubling property and symmetry in $x$ and $y$,
together with~\eqref{eq-restrict} and~\eqref{eq-ass-Y-X-Om}, imply that
\begin{align*}
  [u]_{\theta,p,Y}^p
 & \simeq  [u]_{\theta,p,X}^p
  +     \int_{\Om} u(x)^p \int_{Y \setm X}   I(x,y) \,d\nu(y)\,d\nu(x) \\
  & \simle [u]_{\theta,p,X}^p
  +     \int_{\Om} u(x)^p \int_{X \setm \Om}   I(x,y) \,d\nu(y)\,d\nu(x) 
   \simeq  [u]_{\theta,p,X}^p.
       \qedhere
\end{align*}
\end{proof}

\section{Hyperbolic fillings and capacities on them}
\label{sect-hypfill}

\emph{In this section, we let $Z$ be a compact metric space with
$0 <\diam Z <1$ and equipped with a doubling measure $\nu$.
Let $x_0\in Z$ be fixed.}

\medskip

Hyperbolic fillings will be one of our main tools when obtaining
precise estimates for condenser capacities, based on results 
from Bj\"orn--Bj\"orn--Lehrb\"ack~\cite{BBLeh1} and~\cite{BBLehIntGreen}.
We follow the construction of the hyperbolic filling in
Bj\"orn--Bj\"orn--Shan\-mu\-ga\-lin\-gam~\cite{BBShyptrace} as follows:
Fix two parameters $\alp,\tau>1$
and let $X$ be a hyperbolic filling of $Z$, constructed with these parameters.
More precisely, fix $z_0\in Z$ and set $A_0=\{z_0\}$.
Note that $Z=\BZ(z_0,1)$.
By a recursive construction
using Zorn's lemma or the Hausdorff maximality principle,
for each positive integer $n$ we can choose 
a maximal $\al^{-n}$-separated set $A_n\subset Z$
such that $A_n\subset A_m$ when $m\ge n \ge 0$. 
A set $A\subset Z$ is \emph{$\al^{-n}$-separated} if $d_Z(z,z')\ge
\al^{-n}$ whenever $z,z'\in A$ are distinct.
Then the balls $\BZ(z,\tfrac12\al^{-n})$, $z\in A_n$, are pairwise disjoint.
Since $A_n$ is maximal, the balls $\BZ(z,\al^{-n})$, $z\in A_n$, cover $Z$.

We define the ``vertex set''
\begin{equation*} 
   V=\bigcup_{n=0}^\infty V_n,
   \quad \text{where }
  V_n=\{(z,n): z \in A_n\}.
\end{equation*}
The vertices $v=(x,n)$ and $v'=(y,m)$ form an edge 
(denoted $[v,v']$) in the \emph{hyperbolic filling} $X$ of $Z$
if and only if $|n-m|\le 1$ and
\begin{align*}
\tau \BZ(x,\alpha^{-n})\cap \tau \BZ(y,\alpha^{-m}) \ne\emptyset, 
&\quad \text{if } m=n, \label{eq-tilde-m-n}  \\
\BZ(x,\alpha^{-n})\cap \BZ(y,\alpha^{-m}) \ne\emptyset, 
&\quad \text{if } m=n\pm 1. 
\end{align*}

The hyperbolic filling $X$, seen as a metric space with edges of unit length,
is a Gromov hyperbolic space.
Its uniformization $X_\eps$ with parameter 
$\eps=\log \alp$ is given by the uniformized metric 
\[
d_\eps(x,y) = \inf_\ga \int_\ga e^{-\eps d(\cdot,v_0)}\,ds = \inf_\ga \int_\ga \al^{-d(\cdot,v_0)}\,ds,
\]
where $d(\,\cdot\,,v_0)$ denotes the graph distance to
the root $v_0=(z_0,0)$ of the hyperbolic filling,
$ds$ denotes the arc length,
 and the
infimum is taken over all paths in $X$ joining $x$ to $y$.
We let 
\[
\clXeps=X_\eps \cup \dXeps
\]
 be the completion of $X_\eps$
and equip it with the measure $\mu_\be$ as in~\cite[Section~10]{BBShyptrace},
with 
\[
\be =\eps (1-\theta)p.
\]
Roughly, $\mu_\be$ is obtained by smearing out the measure $\nu(B(x,\al^{-n}))$ 
to the edges adjacent to the vertex $(x,n)\in V$.
Note that $e^\eps=\al$ and that $\sigma$, appearing in various places 
in~\cite{BBShyptrace}, is in our case
\[
    \sigma = \frac{\eps}{\log \alp} = 1.
\]
By \cite[Proposition~4.4]{BBShyptrace},  $Z$ and $\dXeps$
are biLipschitz equivalent (since $\sigma =1$)
and we will therefore identify them as sets.
However, the metrics are different.
More precisely, by \cite[Proposition~4.4]{BBShyptrace},
\begin{equation} \label{eq-Psi-new}
C_1 d_Z(x,y) \le d_\eps(x,y) \le C_2 d_Z(x,y)
  \quad \text{for all $x,y\in Z$},
\end{equation}
where  
$C_1=1/2\tau\alp$, $C_2=4\alp^{(l+1)}/\eps$
and $l$ is the smallest nonnegative integer such that $\al^{-l}\le\tau-1$.

Clearly, $Z$ is uniformly perfect at $x_0$ if
and only if $\dXeps$ is uniformly perfect at $x_0$
(with comparable constants $\ka$ and $\ka_\eps$).
Moreover, if $\Om \subset Z$ is open and $E \Subset \Om$,
then
\begin{equation} \label{eq-cptZ-cptXeps}
\cptZ(E,\Om) \simeq 
  \cptXeps(E,\Om).
\end{equation}
Note however that because of \eqref{eq-Psi-new}, if $E$
and $\Om$ are balls with respect to $Z$, they will not in general
be balls with respect to $\bdy_\eps X$, which needs to be taken
into account when estimating the capacity of annuli.

We will need the Newtonian (Sobolev) space on $\clXeps$ and its
Sobolev and condenser capacities, which we now introduce,
see \cite{BBbook} or \cite{BBShyptrace} for further details.

A property holds for \emph{\p-almost every curve}
in $\clXeps$
if the curve family $\Ga$ for which it fails has zero \p-modulus,
i.e.\ there is $\rho\in L^p(\clXeps)$ such that
$\int_\ga \rho\,ds=\infty$ for every $\ga\in\Ga$.
A measurable
  function $g:\clXeps \to [0,\infty]$ is a \p-weak \emph{upper gradient}
of $u:X \to [-\infty,\infty]$
if for \p-almost all rectifiable curves
$\gamma: [0,l_{\gamma}] \to X$,
\begin{equation*} 
        |u(\gamma(0)) - u(\gamma(l_{\gamma}))| \le \int_{\gamma} g\,ds,
\end{equation*}
where the left-hand side is $\infty$
whenever at least one of the
terms therein is infinite. 
If $u$ has a \p-weak upper gradient in $L^p(\clXeps)$, then
it has a \emph{minimal \p-weak upper gradient}
$g_u \in L^p(\clXeps)$
in the sense that
$g_u \le g$ a.e.\
for every \p-weak upper gradient $g \in L^p(\clXeps)$ of $u$.

For measurable $u:X \to [-\infty,\infty]$, we let
\[
        \|u\|_{\Np(\clXeps)} = \biggl( \int_{\clXeps} |u|^p \, d\mu
                + \inf_g  \int_{\clXeps} g^p \, d\mu \biggr)^{1/p},
\]
where the infimum is taken over all \p-weak upper gradients of $u$.
The \emph{Newtonian space} on $\clXeps$ is
\[
        \Np (\clXeps) = \{u: \|u\|_{\Np(\clXeps)} <\infty \}.
\]
Note that functions in $\Np$ are defined pointwise
everywhere, not only up to a.e.-equivalence classes.

The  \emph{Sobolev capacity} of $E\subset \clXeps$ is
\[
\CpXeps(E)=\inf_{u}\|u\|_{\Np(\clXeps)}^p,
\]
where the infimum is taken over all $u \in \Np(\clXeps)$ such that
$u=1$ on $E$.
The \emph{condenser capacity} of $E \subset \Om$ with respect to
an open set $\Om \subset \clXeps$ is
\[
\cpXeps(E,\Om) = \inf_{u}\int_{\clXeps} g_{u}^p\, d\mu,
\]
where the infimum is taken over all $u \in \Np(\clXeps)$
such that $u=1$ on $E$ and $u=0$ outside $\Om$.
For both capacities we call such $u$ admissible.

By \cite[Theorem~10.3]{BBShyptrace}, $\mube$ is doubling and supports
a $1$-Poincar\'e inequality on $\clXeps$,
i.e.\ there exist $C,\lambda>0$ such that for each ball $B=\BXeps(x,r)$  and for 
all integrable functions $u$ and $1$-weak upper gradients $g$ of $u$ on $\la B$, 
\begin{equation} \label{eq-PI}
\vint_{B}|u-u_{B}|\, d\mube \le C r \vint_{\la B}g\, d\mube,
\end{equation}
where 
\[
u_B:= \vint_B u\,d\mube = \frac{1}{\nu(B)} \int_B u\,d\mube.
\]
As $\clXeps$ is geodesic, the dilation constant in the $1$-Poincar\'e
inequality can be chosen to be $\la=1$
and moreover $\clXeps$ supports  a $(p,p)$-Poincar\'e
inequality (i.e.\ \eqref{eq-PI} with averaged $L^p$-norms on both sides)
with dilation $\la=1$, see e.g.\
\cite[Theorem~4.39]{BBbook}.
It thus follow from
  Bj\"orn--Bj\"orn--Shan\-mu\-ga\-lin\-gam~\cite[Corollary~1.3]{BBS-qcont}
and \cite[Theorems~6.7\,(vii) and~6.19\,(vii)]{BBbook}
that $\CpXeps$ and $\cpXeps$ are outer capacities.

Another consequence of \cite[Theorem~10.3]{BBShyptrace} is that
for every $r \le 2 \diam_\eps \clXeps$ and $x \in Z$,
\begin{equation} \label{eq-mube-nu}
   \mube(\BXeps(x,r))
   \simeq r^{\be/\eps} \nu(\BZ(x,r))
   = r^{(1-\theta)p} \nu(\BZ(x,r)).
\end{equation}

From~\eqref{eq-Psi-new} and~\eqref{eq-mube-nu}
it follows that the exponent sets are the same for
$Z$ and $\dXeps$, and that, for $q>0$,
\[
   q \in \lQo^Z \eqv q + \frac{\be}{\eps}= q+(1-\theta)p \in \lQo^{\clXeps},
\]
and similarly for the other exponent sets.   
Here we consider the exponent sets around $x_0 \in Z$.
Moreover, if $Z$ is uniformly perfect at $x_0$,
then the doubling property implies that  all the exponent sets 
for $\nu$ and $\mu_\be$
are nonempty, see \cite[(2.3)]{BBLeh1}.
Hence
\begin{equation}   \label{eq-transfer-qo-Z-X}
   \lqoZ = \lqoX - \frac{\be}{\eps}= \lqoX - (1-\theta)p,
\end{equation}
and similarly for the other exponents.   
In particular,
\begin{equation} \label{eq-lqoX-Z}
   p<\lqoX \eqv \theta p < \lqoZ.
\end{equation}

We are now ready to estimate
capacities on $\bdyXeps$ in terms of capacities on $\clXeps$,
with the aim to later translate them to capacities on the
original space $Z$.
The comparison constants in this section 
are independent of the choice of $x_0$,
and
depend only on $\theta$, $p$, $C_\nu$, $\alp$ and $\tau$,
  unless said otherwise.

\begin{lem} \label{lem-cpt<cp-E}
Assume that $E \subset \BdXeps_R$.
Then
\begin{equation*} 
  \cptXeps(E,\BdXeps_{2R}) \simle \cpXeps(E,\BXeps_{2R}).
\end{equation*}
\end{lem}

\begin{proof}
As both capacities are outer, we may assume
that $E$ is open in $Z$.
Let $u \in \Np(\clXeps)$ be admissible for $\cpXeps(E,\BXeps_{3R/2})$.
Then the restriction $u|_Z$ is admissible for  $\cptXeps(E,\BdXeps_{2R})$,
and by \cite[Theorem~11.3]{BBShyptrace},
\[
\cptXeps(E,\BdXeps_{2R})
\le [u|_Z]_{\theta,p,\dXeps}^p
\simle \|g_u\|_{L^p(\clXeps)}^p.
\]
Taking infimum over all $u$ admissible for $\cpXeps(E,\BXeps_{3R/2})$
shows that
\[
\cptXeps(E,\BdXeps_{2R}) \simle \cpXeps(E,\BXeps_{3R/2})
\simeq \cpXeps(E,\BXeps_{2R}),
\]
where the last comparison follows from
\cite[Lemma~11.22]{BBbook}.
\end{proof}

The following lemma controls how function values spread from
$Z$ to the hyperbolic filling.
This property will be essential for obtaining a reverse estimate 
to Lemma~\ref{lem-cpt<cp-E}.

\begin{lem}  \label{lem-shadow}
Let $u\in B^\theta_p(Z)$ be such that $u=b$ in $\BdXeps(x,Lr)$,
where $L=1+\alp(1+\eps + C_2\eps)$ with $C_2$ as in~\eqref{eq-Psi-new}.
Let $U$ be the extension  of $u$ to $X_\eps$, given by 
\begin{equation} \label{eq-V-def}
  U((z,n)):=\vint_{\BZ(z,\alpha^{-n})} u\, d\nu,
  \quad \text{if } (z,n)\in V\subset X,
\end{equation}
extended piecewise linearly {\rm(}with respect to $d_\eps${\rm)} 
to each edge in $X_\eps$, and then by
\begin{equation} \label{eq-U-bdyeps}
 U(x):=\limsup_{r \to 0} \vint_{B_{X_\eps}(x,r)} U \, d\mu_\be 
\quad \text{for $x \in \bdy_\eps X$}.
\end{equation}
Then $U \equiv b$ in $\BXeps(x,r)$.
\end{lem}

\begin{proof}
Let $y \in \BXeps(x,r) \setm Z$.
Then $y$ belongs to an edge  $[v_1,v_2]$, where $v_1=(x_1,n_1)$ and $v_2=(x_2,n_2)$
are vertices in the hyperbolic filling. 
We can assume that $n_1 \le n_2\le n_1+1$.
Then for $j=1,2$, since $\al=e^\eps$,
\[
d_\eps(y,v_j) \le \int_{0}^{1} \al^{- n_1} \, dt = \al^{-n_1}
\quad \text{and}   \quad
d_\eps(v_j,x_j) = \frac{\al^{-n_j}}{\eps} \le \frac{\al^{-n_1}}{\eps}.
\]
Since also
\[
r > d_\eps(y,x) \ge \dist_\eps(y,Z)
   \ge \int_{n_2}^\infty \al^{-t} \, dt \ge \frac{\al^{-n_1-1}}{\eps},
\]
we have for all $z\in \BZ(x_j,\al^{-n_j})$, $j=1,2$,
that using also \eqref{eq-Psi-new},
\begin{align*}
  d_\eps(x,z)
  &   < d_\eps(x,y) + d_\eps(y,v_j) + d_\eps(v_j,x_j) + C_2 \alp^{- n_j} \\
  &< r + \al^{-n_1}\biggl(1+\frac{1}{\eps} + C_2\biggr)
  < r + \al \eps r\biggl(1+\frac{1}{\eps} + C_2\biggr)
  = Lr,
 \end{align*}
and thus $u(z)=b$ by assumption.
It follows from~\eqref{eq-V-def} that $U(x_j)=b$, $j=1,2$, and hence also $U(y)=b$.
For $y \in \BXeps(x,r) \cap Z$, the claim follows from~\eqref{eq-U-bdyeps}.
\end{proof}

\begin{thm} \label{thm-cpt-cp-E}
  Assume that $Z$ is uniformly perfect at $x_0$
with constant $\ka$, and that $E \subset \BdXeps_R$.
If $\BdXeps_{3R} \ne Z$
then
\begin{equation} \label{eq-cpt-cp-E}  
\cptXeps(E,\BdXeps_{2R}) \simeq \cpXeps(E,\BXeps_{2R}),
\end{equation}
with comparison constants also depending on $\ka$.
\end{thm}

\begin{proof}
The ``$\simle$'' inequality follows from Lemma~\ref{lem-cpt<cp-E},
so it remains to show the ``$\simge$'' inequality.
As both capacities are outer,  we may assume
that $E$ is open in $Z$.
Let $u$ be admissible for  $\cptXeps(E,\BdXeps_{2R})$.
Consider the extension $U$ to $X_\eps$ given by \eqref{eq-V-def}
and~\eqref{eq-U-bdyeps}.
It then follows from \cite[Theorem~12.1]{BBShyptrace} that
$U=u$ $\nu$-a.e.\ in $\bdy_\eps X$ and
\begin{equation} \label{eq-V}
  \int_{\clXeps} g_{U}^p \, d\mube
  \simle [u]_{\theta,p,\dXeps}^p.
\end{equation}
As $E$ is open, it easily follows (e.g.\ from 
  Lemma~\ref{lem-shadow} and \eqref{eq-U-bdyeps}) that $U \equiv 1$ on $E$.
Moreover $0 \le U \le 1$ on $\clXeps$.

Let next
$\eta : \clXeps \to [0,1]$ be a $2/R$-Lipschitz cut-off function with
$\supp \eta \Subset \BXeps_{2R}$ such that
$\eta=1$ in $\BXeps_R$.
Then, by \cite[Theorem~2.15]{BBbook},
\[
    g_{\eta U} \le \eta g_U + U g_\eta \le g_U + \frac{2U}{R}.
\]
Since $\eta U$ is admissible for $\cpXeps(E,\BXeps_{2R})$, we have
\begin{equation}  
  \cpXeps(E,\BXeps_{2R})
  \le \int_{\BXeps_{2R}} g_{\eta U}^p \, d\mube
  \label{eq-est-with-V-eta}
\simle  \int_{\BXeps_{2R}} g_{ U}^p \, d\mube
    + \frac{1}{R^p} \int_{\BXeps_{2R}} U^p \, d\mube.
\end{equation} 

In view of~\eqref{eq-V}, it therefore suffices to
estimate the last term in \eqref{eq-est-with-V-eta}
using the first integral on the right-hand side.
To this end, let $B=\BXeps_{4\kaeps R}$, where $\kaeps$ 
is the uniform perfectness constant of $\bdyXeps$ at $x_0$.
We will use that
\[
  \frac{\mube(B \setm\spt U)}{\mube(B)}\ge \Theta >0,
\]
where $\Theta$ is independent of $U$ and $B$ and
only depends on $\eps$, $\kaeps$ and $C_\mube$.
We postpone the verification of this to the end of the proof and first
show how it leads us to conclude the proof.
The Minkowski inequality yields
\begin{equation}   \label{eq-est-with-V-B}
   \biggl( \vint_{B} U^p \, d\mube \biggr)^{1/p} 
   \le \biggl( \vint_{B} |U-U_B|^p \, d\mube \biggr)^{1/p} + U_B,
\end{equation}
where
\begin{align*}
U_B &:= \vint_B U \, d\mube 
  = \vint_{B \setm \supp U} |U-U_B|  \, d\mube \\
  & \le \frac{\mube(B)}{\mube(B \setm \supp U)}
      \vint_{B} |U-U_B|  \, d\mube
  \le  \frac{1}{\Theta}
 \biggl( \vint_{B} |U-U_B|^{p}  \, d\mube \biggr)^{1/p}.
\end{align*}
Inserting this into \eqref{eq-est-with-V-B} and using the $(p,p)$-Poincar\'e
inequality for $\mube$ gives
\[
\int_{B} U^p \, d\mube 
\simle \vint_{B} |U-U_B|^{p}  \, d\mube
\simle R^p \int_B g_{ U}^p \, d\mube.
\]
Together with~\eqref{eq-V} and \eqref{eq-est-with-V-eta} 
the last estimate implies that
\begin{equation} \label{eq-3}
  \cpXeps(E,\BXeps_{2R})
   \simle \int_{\BXeps_{2R}} g_U^p \, d\mube  
    +  \int_{B} g_U^p \, d\mube   
  \simle [u]_{\theta,p,\dXeps}^p.
\end{equation} 
Taking infimum over all $u$ admissible for 
$\cptXeps(E,\BdXeps_{2R})$ shows the
  ``$\simge$'' inequality in \eqref{eq-cpt-cp-E}.

It remains to show that $\Theta > 0$.
By the uniform perfectness and the fact that
$\BdXeps_{3 R} \ne Z$, there is some
$x \in \BdXeps_{3 \kaeps R} \setm \BdXeps_{3 R}$.
Then $u=0$ in $\BdXeps(x,R)$ and hence by 
Lemma~\ref{lem-shadow},
$U=0$ in $\BXeps(x,R/L)$. 
From this and the doubling property of $\mube$ we see that
\[
 \frac{\mube(B \setm\spt U)}{\mube(B)} 
\ge  \frac{\mube(\BXeps(x,R/L))}{\mube(B)} 
\ge \Theta> 0,
\]
where $\Theta$ only depends on $\eps$, $\kaeps$ and $C_\mube$.
\end{proof}  

Since we will be interested in the Besov capacity of annuli in $Z$, 
we next relate it to the capacity of annuli in $\clXeps$.

\begin{thm} \label{thm-cpt-cp-B}
  Assume that $Z$ is uniformly perfect at $x_0$
 with constant $\ka$.
  Let $0 < 2r \le R$ and 
  $L=\alp(1+\eps + C_2\eps)$ 
    as in Lemma~\ref{lem-shadow}.
  Assume that $\BdXeps_{3R/2} \ne Z$.
Then 
\begin{equation} \label{eq-cpt-cp-B}  
 \cptXeps(\BdXeps_r,\BdXeps_R)
 \simge  
\cpXeps(\BXeps_{r/L},\BXeps_R),
\end{equation}
with comparison constant also depending on $\ka$.
\end{thm}

\begin{proof}
We proceed
as in the proof of Theorem~\ref{thm-cpt-cp-E} with
$E=\BdXeps_r$ and $2R$ replaced by $R$.
Lemma~\ref{lem-shadow} shows that the function
$U$ constructed in \eqref{eq-V-def}
and~\eqref{eq-U-bdyeps} satisfies
$U \equiv 1$ in $\BXeps_{r/L}$ and is thus admissible
for $\cpXeps(\BXeps_{r/L},\BXeps_R)$,
i.e.\  we can replace $E$ by $\BXeps_{r/L}$
in~\eqref{eq-3}.
Taking infimum over all $u$ admissible for 
$\cptXeps(\BdXeps_r,\BdXeps_R)$ shows \eqref{eq-cpt-cp-B}. 
\end{proof}  

\section{Enlarging \texorpdfstring{$Y$}{Y}}
\label{sect-ext-Y}

\emph{In this section we assume that $Y$ is
a compact metric space, equipped
with a doubling measure~$\nu$, and let $x_0\in Y$ be fixed.}

\medskip

Our aim is to embed $Y$ into a  suitable larger metric space $Z$.
We will do this recursively, but in this section we only do the first step.

As $Y$ is compact there is a point $x_1$
such that $d(x_1,x_0)=\max_{x \in Y}  d(x,x_0)$.
Let $Y'=(Y',d',\nu')$ be a copy of $Y=(Y,d,\nu)$,
where we identify $x_1$ with its copy, but do not identify any other points.
Equip  $\Yhat=Y \cup Y'$ with the measure 
\[
\nuhat(A)=\nu(A\cap Y)+\nu'(A\cap Y')
\]
and the metric $\dhat$ so that
\[
  \dhat(x,y)=\begin{cases}
    d(x,x_1)+d'(y,x_1), & \text{if }
    x \in Y \text{ and } y \in Y', \\
    d(x,y) & \text{if }     x,y \in Y, \\ 
    d'(x,y) & \text{if }     x,y \in Y'.
  \end{cases}
\]

\begin{lem} \label{lem-nu-Y-Yhat}
The measure $\nuhat$ is doubling on $\Yhat$
with doubling constant $C_\nuhat \le 2 C_\nu$
and satisfies
\begin{equation} \label{eq-nuhat}
\nu(\BYY(x,r))
\le  \nuhat(\BYhat(x,r))
\le  2\nu(\BYY(x,r))
\quad \text{if $x \in Y$ and $r>0$}.
\end{equation}

Moreover, if $Y$ is uniformly perfect at $x_0$ with constant $\ka$,
then $\Yhat$ is uniformly perfect at $x_0$ with constant
$\kat=\max\{\ka,2\}$.
\end{lem}

\begin{proof}
That \eqref{eq-nuhat} holds follows directly from the construction.
  A similar formula holds if $x \in Y'$.
It follows that $\nuhat$ is doubling with $C_{\nuhat} \le 2 C_\nu$.

As for the uniform perfectness, let $r>0$ be such that
$\BYhat_{\kat r} \ne \Yhat$.
Then $\kat r \le 3 d(x_0,x_1)$ and hence $r \le \tfrac32 d(x_0,x_1)$.
If $r \le d(x_0,x_1)$ then $x_1 \in Y \setm \BYY_{r}$ and thus
there is
\[
y \in \BYY_{\ka r} \setm \BYY_r \subset \BYhat_{\kat r} \setm \BYhat_r.
\]
On the other hand, if $d(x_0,x_1) < r \le \tfrac32 d(x_0,x_1)$,
then $\BYhat_{\kat r} \setm \BYhat_r$ contains the copy of $x_0$ in $Y'$. 
\end{proof}  

The constant $2$ in $\kat$ in Lemma~\ref{lem-nu-Y-Yhat} is optimal
as seen by the following example:
Let $Y=[-1,0] \cup \{1\}$ with
$x_0=0$ and $x_1=1$. In this case
$Y$ is uniformly perfect at $0$ with any constant $\ka >1$,
but $\Yhat$ is only uniformly perfect at $0$ with constant $\kat \ge 2$.

From now on we call the distance $d$ and the measure $\nu$ also on $\Yhat$.

\begin{lem} \label{lem-cp-restrict-Yhat}
 Let $\Om \subset B_{d(x_0,x_1)/2}$ be open
  and $E \Subset \Om$.
Then 
\[
\cptY(E,\Om) \simeq \cptYhat(E,\Om),
\]
with comparison constants
depending only on $\theta$, $p$  and $C_\nu$.
\end{lem}

\begin{proof}
Lemma~\ref{lem-nu-Y-Yhat} shows that \eqref{eq-restrict} in
Lemma~\ref{lem-cp-restrict-new} holds for the spaces $Y\subset \Yhat$.
By the doubling property of $\nu$,
\[
\nu(\Yhat \setm Y) \simeq \nu(B(x_1, \tfrac12 d(x_0,x_1)) \simeq \nu(Y \setm \Om)
\]
Since for all $x\in\Om$, $y\in \Yhat \setm Y$
and $y'\in Y\setm\Om$,
\[
d(x,y) \simeq d(x_0,x_1) 
\quad \text{and} \quad d(x,y') \le \tfrac32 d(x_0,x_1),
\]
the statement follows
from Lemma~\ref{lem-cp-restrict-new}.
\end{proof}  

\section{From unbounded to bounded spaces}
\label{sect-unbdd}

\emph{In this section, we let $Y$ be a metric space equipped with a doubling measure
$\nu$ and fix $x_0 \in Y$.}

\begin{lem}  \label{lem-doubl-begr}
Let $Y_0=\{x_0\}$ and $\de>0$.
For $n=0,1,\ldots$\,, let 
\[
Y_{n+1} =  \bigcup_{x\in Y_n} B(x, 2^{-n}\de)
\quad \text{and}   \quad
Y'= \overline{\bigcup_{n=0}^\infty Y_n}.
\]
Also let $\nu':= \nu|_{Y'}$.

Then the following hold\/\textup{:}
\begin{enumerate}
\item \label{k1}
  $\BYY_{\de} \subset Y' \subset \overline{\BYY_{2\de}}$,
\item \label{k2}
$\nu'$ is doubling with $C_{\nu'} \le C_\nu^6$,  
\item \label{k3}
for all $x\in Y'$ and $0 < r < 2\diam Y'$,
\[
\nu'(\BYp(x,r)) \le \nu(\BYY(x,r)) \le C_d^5 \nu'(\BYp(x,r)),
\]
\item \label{k4}
  if $Y$ is uniformly perfect at $x_0$ with constant $\ka$,
  then $Y'$ is uniformly perfect at $x_0$ with constant $\ka'=\max\{\ka,2\}$.
\end{enumerate}
\end{lem}

\begin{proof}
That \ref{k1} holds is clear from the construction.

\ref{k3}
The first inequality is obvious.
By \ref{k1}, $r < 2 \diam Y' \le 8 \de$.
Let $r'= 2^{3-k} \de$, where $k \ge 0$ is an integer such that
$ \tfrac12 r < r' \le  r$.
Then there is $x' \in Y_k$ such that $d(x,x')< 2^{1-k} \de = \tfrac14 r'$.
Thus $\BYY(x',\tfrac18 r') \subset Y_{k+1} \subset Y'$ and so
\begin{align*}
\nu(\BYY(x,r))
& \le \nu(\BYY(x',4r'))
\le C_\nu^5 \nu(\BYY(x',\tfrac18 r')) \\
& =  C_\nu^5 \nu'(\BYp(x',\tfrac18 r'))
\le  C_\nu^5 \nu'(\BYp(x,r)).
\end{align*}

\ref{k2}
Let $x \in Y'$ and $r >0$.
If $r < 2 \diam Y'$, then by \ref{k3},
\[
\nu'(\BYp(x,2r))
\le \nu(\BYY(x,2r))
\le C_\nu \nu(\BYY(x,r))
\le C_\nu^6 \nu'(\BYp(x,r)).
\]
If instead $r \ge 2 \diam Y'$, then with $\rt=\diam Y'$,
\[
\nu'(\BYp(x,2r))
= \nu'(\BYp(x,2\rt))
\le C_\nu^6 \nu'(\BYp(x,\rt))
\le C_\nu^6 \nu'(\BYp(x,r)).
\]

\ref{k4}
Let $r>0$ be such that $\BYp_{\ka'r} \ne Y'$.
Then 
$Y \setm \BYY_{\ka' r}     \supset Y' \setm \BYp_{\ka'r}\ne\emptyset$
  and $\ka' r \le 2 \de$.
Hence, if $\ka' r \le \de$ then
there is $z \in \BYY_{\ka' r} \setm \BYY_{r} = \BYp_{\ka' r} \setm \BYp_{r}$.
So we may assume that $\de < \ka' r \le 2 \de$.
As $\BYp_{\ka'r} \ne Y' \supset Y_1=\BYp_{\de}$ we
see that $Y_2 \setm Y_1 \ne \emptyset$.
Therefore there are
$x_1 \in Y_1$ and $x_2 \in Y_2 \setm Y_1$ with $d(x_1,x_2) < \tfrac12 \de$.

Assume for a contradiction that $\BYp_{\ka' r} \setm \BYp_r = \emptyset$.
Since $ r \le \de < \ka' r$ we must have
$d(x_1,x_0) < r$ and hence also
\[
d(x_2,x_0) \ge \ka' r \ge 2r > 2 d(x_1,x_0).
\]
Thus,
\[
\tfrac12 \de > d(x_1,x_2) \ge  d(x_2,x_0) - d(x_1,x_0)
> \tfrac12 d(x_2,x_0) \ge \tfrac12\ka' r >\tfrac12 \de,
\]
a contradiction.
Hence $\BYp_{\ka' r} \setm \BYp_r \ne \emptyset$.
\end{proof}

The following lemma shows that for the condenser capacity,
the (possibly unbounded) space $Y$ can be effectively replaced by 
the bounded space $Y'$.

\begin{lem} \label{lem-cpt-restrict-ball}
Let $Y'$ be the space constructed in Lemma~\ref{lem-doubl-begr} with
parameter $\de>0$.
Assume that $Y$ is uniformly perfect at $x_0$
  with constant $\ka$.
  Let $R=\de/2 \ka$,
  $\Om \subset \BYY_R$ be open and  $E \Subset \Om$.
Then 
\[
\cptYprime(E,\Om) \simeq \cptY(E,\Om),
\]
with comparison constants
depending only on $\theta$, $p$, $C_\nu$ and $\ka$.
\end{lem}  

The assumption of uniform perfectness cannot be dropped
since  
$\cptYprime(E,\Om)=0$  if $Y' \setm\Om=\emptyset$.

\begin{proof}
We shall use Lemma~\ref{lem-cp-restrict-new}.
If $Y'=Y$, there is nothing to prove, so assume that $Y' \ne Y$.
Then $\BYY_{2 \ka R} =\BYY_{\de} = Y_1 \ne Y$.
By the uniform perfectness of $Y$, there is some
$z \in \BYY_{2 \ka R} \setm \BYY_{2 R} \subset Y_1$.
Then $\BYY(z,R) \subset Y_2 \setm \Om \subset Y' \setm \Om$.
Let $x \in \Om$ and $y\in Y'$.
Then 
\[
d(x,y) \le 2\de + R = (4\ka+1)R,
\]
and hence, using that
$\nu'=\nu|_{Y'}$ 
is doubling by Lemma~\ref{lem-doubl-begr}, we obtain
\begin{align*}
\nu'(\BYp(x,d(x,y)))  &\simle \nu'(\BYp(x,R))  \\
&\le \nu'(\BYp(z,2(\ka+1)R)) 
\simeq \nu'(\BYp(z,R))
= \nu(\BYY(z,R)).
\end{align*}
Thus, with $I(x,y)$ as in Lemma~\ref{lem-cp-restrict-new},
\[ 
\int_{Y' \setm \Om} I(x,y) \,d\nu(y)
\simge \int_{\BYY(z,R)} \frac{d\nu(y)}{R^{\theta p}\nu(\BYY(z,R))}
= R^{-\theta p}.
\] 
On the other hand, for
$y\in A^j:=\BYY_{2^{j+1}\de} \setm \BYY_{2^j\de}$, $j=0,1,\ldots$\,,
we have
\[
d(x,y) \simeq 2^jR \quad \text{and} \quad
\nu(\BYY(x,d(x,y))) \simge \nu(A_j).
\]
Hence
\begin{align*}
  \int_{Y \setm Y'}
  I(x,y) \,d\nu(y)
   & \le  \sum_{j=0}^\infty \int_{A^j}  I(x,y) \,d\nu(y) \\
   &\simle \sum_{j=0}^\infty \frac{1}{(2^{j}R)^{\theta p}}
   \simeq {R^{-\theta p}} \simle \int_{Y' \setm \Om} I(x,y) \,d\nu(y).
\end{align*}
An application of Lemma~\ref{lem-cp-restrict-new},
together with Lemma~\ref{lem-doubl-begr}\,\ref{k3}, 
concludes the proof.
\end{proof}

\section{Proof of Theorem~\ref{thm-metaest}}
\label{sect-sharp-est}

\begin{lem}  \label{lem-comp-int}
Let $0<\Theta_1<\Theta_2< \infty$ and $r>0$.
If $\nu$ is doubling, then
\[
\int_{\Theta_1r}^{\Theta_2r} \biggl( \frac{\rho^{\theta p}}{\nu(\BY_\rho)} \biggr)^{1/(p-1)} 
\,\frac{d\rho}{\rho}
   \simeq \biggl(  \frac{r^{\theta p}}{\nu(\BY_r)} \biggr)^{1/(p-1)}
\]
with comparison constants depending 
only on $\Theta_1$, $\Theta_2$, $\theta$, $p$ and $C_\nu$.
\end{lem}

\begin{proof}
By the doubling property of $\nu$, we have
$\nu(\BY_\rho)\simeq \nu(\BY_r)$ for all
$\Theta_1 r\le \rho\le\Theta_2r$.
The statement now follows by direct calculation of the integral.
\end{proof}

\begin{remark} \label{rmk-thm-1.1}
The comparison constants in Theorem~\ref{thm-metaest}
are independent of the choice of $x_0$.
They depend only on $\theta$, $p$, $C_\nu$
and the uniform perfectness constant~$\ka$.
In the proof below, the
constants $C_1$ and $C_2$ (and thus the ultimate
comparison constants) depend on $\alp$ and $\tau$.
To avoid this dependence 
in Theorem~\ref{thm-metaest}, we can e.g.\ let $\alp=\tau=2$.
We have chosen not to fix $\alp$ and $\tau$, so as to show
that our proof is not dependent on fixing them.
\end{remark}

\begin{proof}[Proof of Theorem~\ref{thm-metaest}]
Let $0<C_1<1<C_2$ be the constants appearing
in \eqref{eq-Psi-new}, which only depend on $\alp$,
$\tau$ and $\eps=\log \alp$.
We can assume that $\ka\ge2$.
To be able to use the hyperbolic filling and 
the capacity results from
Section~\ref{sect-hypfill}, we need to use the results
from either Section~\ref{sect-ext-Y} or Section~\ref{sect-unbdd},
depending on if $Y$ is bounded or not.

If $Y$ is bounded, we
use the construction from Section~\ref{sect-ext-Y}
recursively  $N$ times (with $N$ only depending on $C_2/C_1$)
and replace $Y$ by its suitable enlargement $Z$
so that $\BZ_{5C_2R/C_1}\ne Z$. 
Note that the doubling constant of $\nu$ is only enlarged by
a factor depending only on  $N$.

If on the other hand $Y$ is unbounded, we 
   let $Z=Y'$, where $Y'$ is as 
in  Lemma~\ref{lem-doubl-begr} with   $\de=5\ka C_2 R/C_1$. 
  We will also denote the restricted measure by $\nu$.
By Lemma~\ref{lem-doubl-begr},
the doubling constant of $\nu$
is 
in this case only enlarged by the power $6$.
Note that $\BZ_{5 C_2 R/C_1} \subsetneq Z$ by the uniform perfectness condition.
  
The uniform perfectness constant $\ka\ge2$ remains unchanged
in both cases.
Since the left- and right-hand sides in \eqref{eq-metaest} and~\eqref{eq-metaest-0}
scale in the same way,
we may without loss of generality assume that $0 < \diam Z < 1$.

If $Y$ is bounded, we apply Lemma~\ref{lem-cp-restrict-Yhat}
with $E=B_r$ and $\Om=B_R$
several times to the consecutive enlargements of $Y$.
If $Y$ is unbounded, we instead use Lemma~\ref{lem-cpt-restrict-ball}.
In both cases we obtain that
\begin{equation}    \label{eq-capY-sim-capZ}
\cptY(B_r,B_R) \simeq \cptZ(\BZ_r,\BZ_R),
\end{equation}
so it suffices to estimate the latter capacity.
We consider two cases.

If $2C_2r\ge C_1R$, then Proposition~\ref{prop-cpY-Theta}
yields
\[
\cptZ(\BZ_r,\BZ_R) \simeq \frac{\nu(\BY_R)}{R^{\theta p}},
\]
which by  Lemma~\ref{lem-comp-int} is comparable to the
integral in \eqref{eq-metaest}.

If $2C_2r\le C_1R$, we follow Section~\ref{sect-hypfill} and construct a 
hyperbolic filling $X$ of $Z$ with parameters $\alp$ and $\tau$, which we uniformize with
parameter $\eps=\log\al$ and equip with the measure $\mube$, with
$\be=\eps(1-\theta)p$, as in  Section~\ref{sect-hypfill}.
As $\BZ_{5C_2R/C_1}  \ne Z$ we see that $\BdXeps_{5C_2R}  \ne Z$
and thus $\diam \clXeps \ge \diam \dXeps \ge 5C_2R$.
We can then use   Theorem~\ref{thm-cpt-cp-E},
together with \eqref{eq-Psi-new}, \eqref{eq-cptZ-cptXeps}
and \cite[Lemma~11.22]{BBbook}, 
to conclude that
\begin{align}  
\cptZ(\BZ_r,\BZ_R)
& \simle \cptXeps(\BdXeps_{C_2r},\BdXeps_{C_1R})   \nonumber \\
& \simle \cpXeps(\BXeps_{C_2r},\BXeps_{C_1R})
\simeq \cpXeps(\BXeps_{C_2r},\BXeps_{C_2R}).
\label{eq-est-cpt-le-cp} 
\end{align}
Similarly, from Theorem~\ref{thm-cpt-cp-B}, \eqref{eq-Psi-new},
\eqref{eq-cptZ-cptXeps} and \cite[Lemma~11.22]{BBbook}
we get
\begin{align}
  \cptZ(\BZ_r,\BZ_R)
  & \simge  \cptXeps(\BdXeps_{C_1r},\BdXeps_{C_2R})   \nonumber \\
 &\simge \cpXeps(\BXeps_{C_1r/L},\BXeps_{C_2R})
\simeq \cpXeps(\BXeps_{C_1r/L},\BXeps_{C_1R/L}),
\label{eq-est-cpt-ge-cp}
\end{align}
where $L$ is as in Theorem~\ref{thm-cpt-cp-B}.

Next,
the comparison~\eqref{eq-mube-nu} between $\mu_\be$ and $\nu$ gives
\[
\biggl( \frac{\rho^{p}}{\mu_\be(\BXeps_\rho)} \biggr)^{1/(p-1)} 
\simeq \biggl( \frac{\rho^{\theta p}}{\nu(\BY_\rho)} \biggr)^{1/(p-1)}.
\]
Together with Bj\"orn--Bj\"orn--Lehrb\"ack~\cite[Theorem~4.2]{BBLehIntGreen} and
the doubling property of $\mu_\be$
this shows that
\begin{align*}
\cpXeps(\BXeps_{C_2r},\BXeps_{C_2R})
& \simeq \biggl(
\int_{C_2r}^{C_2R} \biggl( \frac{\rho^{p}}{\mu_\be(\BXeps_{\rho})} \biggr)^{1/(p-1)} 
           \,\frac{d\rho}{\rho} \biggr)^{1-p} \\
& \simeq \biggl(
\int_{r}^{R} \biggl( \frac{\rho^{\theta p}}{\nu(\BY_\rho)} \biggr)^{1/(p-1)} 
           \,\frac{d\rho}{\rho} \biggr)^{1-p}.
\end{align*}
Similarly,
\[ 
\cpXeps(\BXeps_{C_1r/L},\BXeps_{C_1R/L})
 \simeq \biggl(
\int_{r}^{R} \biggl( \frac{\rho^{\theta p}}{\nu(\BY_\rho)} \biggr)^{1/(p-1)} 
           \,\frac{d\rho}{\rho} \biggr)^{1-p},
\] 
which together with \eqref{eq-capY-sim-capZ}--\eqref{eq-est-cpt-ge-cp}
concludes the proof of \eqref{eq-metaest}.
The estimate for $\cptY(\{x_0\},\BY_R)$ follows immediately by letting
$r\to0$ in \eqref{eq-metaest}\, since $\cptY$ is an outer capacity.
\end{proof}

\section{Proofs of Theorems~\ref{thm-lqo} and~\ref{thm-S}}

\label{sect-borderline}

\begin{proof}[Proof of Theorem~\ref{thm-lqo}]
  The upper bounds follow directly from Proposition~\ref{prop-cpt-le-new}.
For the lower bounds we first construct $Z$ as in
    the proof of Theorem~\ref{thm-metaest}.
Since the left- and right-hand sides in \eqref{eq-it-1} and~\eqref{eq-it-2}
scale in the same way,
we may without loss of generality assume that $0 < \diam Z < 1$.

    As in~\eqref{eq-est-cpt-ge-cp}, we see that
\begin{equation} \label{eq-cptY-cpX}  
 \cptZ(\BZ_r,\BZ_R)
    \simge  \cptXeps(\BdXeps_{C_1r},\BdXeps_{C_2R})   
 \simge \cpXeps(\BXeps_{C_1r/L},\BXeps_{C_2R}),
\end{equation}
where $L$ is as in Theorem~\ref{thm-cpt-cp-B}.
In \ref{it-1},
it follows from \eqref{eq-lqoX-Z} that $p < \lqoX$.
Hence, by 
Bj\"orn--Bj\"orn--Lehrb\"ack~\cite[Theorem~1.1]{BBLeh1},
\eqref{eq-mube-nu} and the doubling property,
\begin{equation} \label{eq-cpt-cp-Br}  
 \cpXeps(\BXeps_{C_1r/L},\BXeps_{C_2R})
  \simge \frac{\mube(\BXeps_{C_1r/L})}{(C_1r/L)^p}
  \simeq \frac{\nu (B_r)r^{(1-\theta)p}}{r^p}
   = \frac{\nu (B_r)}{r^{\theta p}}.
\end{equation}
In \ref{it-2}, we instead have $p > \uqoX$  and \cite[Theorem~1.1]{BBLeh1},
together with~\eqref{eq-mube-nu} and the doubling property, yields
\begin{equation} \label{eq-cpt-cp-BR}  
 \cpXeps(\BXeps_{C_1r/L},\BXeps_{C_2R}) 
   \simge \frac{\mube(\BXeps_{C_2R})}{(C_2R)^p}
  \simeq   \frac{\nu (B_R)}{R^{\theta p}}.
\end{equation}
Inserting \eqref{eq-cpt-cp-Br} and \eqref{eq-cpt-cp-BR} into~\eqref{eq-cptY-cpX}
and using \eqref{eq-capY-sim-capZ}
proves the lower bounds in~\eqref{eq-it-1} and \eqref{eq-it-2}.

It remains to discuss the sharpness.
Let $0 < 2r < R \le \tfrac14\diam Y$.
If the lower bound in \eqref{eq-it-1} holds, then by 
Proposition~\ref{prop-cpt-le-new},
\[
\frac{\nu (B_r)}{r^{\theta p}} \simle  \cptY(B_r,B_{R})
\simle \frac{\nu (B_R)}{R^{\theta p}},
\]
which immediately implies that $\theta p \in  \lQoY$.
The argument for \eqref{eq-it-2} is similar, using
the upper bound $\nu (B_r)/r^{\theta p}$ from 
Proposition~\ref{prop-cpt-le-new}.

Finally, if $p>1$ then Theorem~\ref{thm-cpt-cp-E},
together with
\eqref{eq-Psi-new}, \eqref{eq-cptZ-cptXeps}, 
\eqref{eq-capY-sim-capZ},
\eqref{eq-it-2} and \eqref{eq-mube-nu}, 
yields
\begin{align*}
\cpXeps(\BXeps_{r},\BXeps_{R}) 
&\simge \cptXeps(\BdXeps_{r},\BdXeps_{R}) 
\simge \cptZ(\BZ_{r/C_2},\BZ_{R/C_1}) \\
&\simeq \cptY(B_{r/C_2},B_{R/C_1})
\simge \frac{\nu(B_R)}{R^{\theta p}} 
= \frac{\nu(\BZ_R)}{R^{\theta p}}
\simeq \frac{\mube(\BXeps_R)}{R^{p}}.
\end{align*}
Theorem~1.3 in Bj\"orn--Bj\"orn--Christensen~\cite{BBChr}, applied to $\clXeps$,
 then implies that 
$p > \uqoX$, which is equivalent to $\theta p > \uqoY$.
\end{proof}

In the borderline cases we have the following result corresponding to 
  Theorem~\ref{thm-lqo}.

\begin{thm} \label{thm-borderline}
Assume that $Y$ is a complete metric space which 
is uniformly perfect at $x_0$ 
and equipped with a doubling 
measure $\nu$. 
Let $p>1$, $0<\theta<1$ and $0 < R_0 \le \tfrac{1}{4} \diam Y$, 
with $R_0$ finite.

Then the following hold for   $0<2r \le R \le R_0$, with 
comparison constants depending on $R_0$,
 but independent of $x_0$, $r$ and $R$.
\begin{enumerate}
\item  \label{bd-1}
If $\theta p =\max\lQoY$, then
\begin{equation} \label{eq-borderline-1}
\frac{\nu(B_r)}{r^{\theta p}} \biggl( \log\frac{R}{r} \biggr)^{1-p}
\simle \cptY(B_r,B_{R})
\simle \frac{\nu(B_R)}{R^{\theta p}} \biggl( \log\frac{R}{r} \biggr)^{1-p}.
\end{equation}
\item  \label{bd-2}
If $\theta p = \min\uQoY$, then
\begin{equation} \label{eq-borderline-2}
\frac{\nu(B_R)}{R^{\theta p}} \biggl( \log\frac{R}{r} \biggr)^{1-p}
\simle \cptY(B_r,B_{R})
\simle \frac{\nu(B_r)}{r^{\theta p}} \biggl( \log\frac{R}{r} \biggr)^{1-p}.
\end{equation}
\end{enumerate}
\bigskip
%

Moreover, if the lower bounds in \eqref{eq-borderline-1} and \eqref{eq-borderline-2}
hold, then $\theta p\le\sup\lQoY$ and $\theta p\ge\inf\uQoY$, respectively.
\end{thm}

\begin{proof} 
The estimate \eqref{eq-borderline-1} 
follows directly from Theorem~\ref{thm-metaest}
since
\[
\frac{R^{\theta p}}{\nu(B_R)}
\simle \frac{\rho^{\theta p}}{\nu(B_\rho)}
\simle \frac{r^{\theta p}}{\nu(B_r)}
\]
as $\theta p = \max\lQoY$. 
The estimate \eqref{eq-borderline-2}
is shown similarly.

For the last statement, the lower bound in \eqref{eq-borderline-1} 
and Proposition~\ref{prop-cpt-le-new}
imply for all $\eps>0$ that 
\[
\frac{\nu(B_r)}{\nu(B_R)} 
\simle \frac{r^{\theta p} \cptY(B_r,B_R)}{\nu(B_R)} 
\biggl( \log \frac{R}{r} \biggr)^{p-1} 
\simle \Bigl( \frac{r}{R} \Bigr)^{\theta p} 
\biggl( \log \frac{R}{r} \biggr)^{p-1} 
\simle \Bigl( \frac{r}{R} \Bigr)^{\theta p-\eps},
\]
where the implicit constant in the last "$\simle$" depends on $\eps$.
Thus $\theta p-\eps\in \lQoY$ for every $\eps>0$, showing that 
$\theta p\le\sup\lQoY$.
The implication \eqref{eq-borderline-2} $\imp$ $\theta p\ge\inf\uQoY$
is proved similarly.
\end{proof}

\begin{remark} \label{rmk-lQ}
If $Y$ is unbounded,  then 
Theorems~\ref{thm-lqo} and~\ref{thm-borderline}
hold with $R_0=\infty$ if $\lQo$, $\lqo$, $\uQo$  and $\uqo$
are replaced by 
\begin{alignat*}{2}
\lQ  &=\biggl\{q>0 : 
        \frac{\mu(B_r)}{\mu(B_R)}  \simle \Bigl(\frac{r}{R}\Bigr)^q 
        \text{ for } 0 < r < R < \infty
        \biggr\},
        &  \quad\lqq &= \sup \lQ, \\
  \uQ
       &=\biggl\{q>0 : 
       \frac{\mu(B_r)}{\mu(B_R)} 
       \simge  \Bigl(\frac{r}{R}\Bigr)^q 
       \text{ for } 0 < r < R < \infty,
       \biggr\}, 
  & \quad \uq &= \inf \uQ.       
\end{alignat*}
\end{remark}  

\begin{remark} \label{rmk-dep-const-lqo}
The comparison constants in Theorems~\ref{thm-lqo} and~\ref{thm-borderline}
are independent of the choice of $x_0$,
but depend  on $\theta$, $p$, $C_\nu$, $R_0$
and the uniform perfectness constant~$\ka$.

In Theorem~\ref{thm-lqo}\,\ref{it-1} 
they also depend on the choice of $q\in(\theta p, \lqo)$ 
from the proof of \cite[Proposition~6.1]{BBLeh1} leading to the
estimate~\eqref{eq-cpt-cp-Br},
and on the comparison constant appearing in the definition of $q \in \lQo$.

Similarly, in Theorem~\ref{thm-lqo}\,\ref{it-2} 
the constants also depend on the choice of $q\in(\uqo,\theta p)$ 
from the proof of \cite[Proposition~6.1]{BBLeh1} leading to the
estimate~\eqref{eq-cpt-cp-BR},
and on the comparison constant appearing in the definition of $q \in \uQo$.

In Theorem~\ref{thm-borderline} the dependence is similar
but with $q=\theta p$.
In Remark~\ref{rmk-lQ}, the dependence is instead in terms of
$\lQ$ and $\uQ$.
\end{remark}

\begin{proof}[Proof of Theorem~\ref{thm-S}]
Let $Z=Y'$, where $Y'$ is as 
in Lemma~\ref{lem-doubl-begr} with $\de=\tfrac15$.
Then $0 < \diam Z < 1$.
(If $Y$ is bounded we may instead let $Z$ be a rescaled version of $Y$.)
Then let $X_\eps$ be the uniformized hyperbolic 
filling for $Z$ constructed in Section~\ref{sect-hypfill}.
By Corollary~\ref{cor-singleton}, it suffices to prove the
statements \ref{it-1-s} and \ref{it-2-s}
for $\CptY(\{x_0\})$, which in turn is comparable to
$\CpXeps(\{x_0\})$ by \cite[Proposition~13.2]{BBShyptrace}.

As in~\eqref{eq-transfer-qo-Z-X}, it follows that 
$p > \inf \uSo^{\clXeps}$
in \ref{it-1-s}, while $p \notin \uSo^{\clXeps}$
or $1<p \in \lSo^{\clXeps}$ in \ref{it-2-s}.
Hence, Proposition~8.2 in \cite{BBLeh1} implies that $\CpXeps(\{x_0\})>0$
in \ref{it-1-s}, and $\CpXeps(\{x_0\})=0$ in \ref{it-2-s}.

When $p>1$, the conclusions can also be derived from
Theorem~\ref{thm-metaest}.
\end{proof}

\end{document}